\documentclass[11pt,reqno]{amsart}

\usepackage[margin=1.5in]{geometry}
\usepackage[T1]{fontenc}

\usepackage{amssymb,amsmath,amsthm,amsfonts,mathrsfs,bm,esint,mathtools}
\allowdisplaybreaks[1]
\numberwithin{equation}{section}

\usepackage{tikz}
\usepackage{float}
\usepackage[font=footnotesize]{caption}

\usepackage{times,comment,soul,xcolor}

\usepackage[numbers,sort&compress]{natbib}
\usepackage[hyperfootnotes=false,pagebackref,pdfstartview=FitV]{hyperref}
\hypersetup{colorlinks,citecolor=blue,linkcolor=red,urlcolor=blue}

\theoremstyle{plain}
\newtheorem{theorem}{Theorem}[section]

\newtheorem{corollary}[theorem]{Corollary}
\newtheorem{lemma}[theorem]{Lemma}

\theoremstyle{definition}

\theoremstyle{remark}
\newtheorem{remark}[theorem]{Remark}

\newcommand{\Z}{\mathbb{Z}}

\newcommand{\T}{\mathbb{T}}
\newcommand{\cH}{\mathcal{H}}
\newcommand{\cL}{\mathcal{L}}
\newcommand{\e}{\varepsilon}
\newcommand{\eps}{\varepsilon}
\DeclareMathOperator{\spn}{span}

\DeclareMathOperator{\diver}{div}

\newcommand{\Norm}[2]{\lVert #1\rVert_{#2}}
\DeclarePairedDelimiter{\norm}{\lVert}{\rVert}
\newcommand{\gd}{\nabla}

\newcommand{\pt}{\partial}
\newcommand{\Om}{\Omega}
\newcommand{\R}{\mathbb R}
\newcommand{\dd}{\,\mathrm{d}}
\newcommand{\dist}{\operatorname{dist}}
\newcommand{\Lop}{\mathcal L}
\newcommand{\AN}{\mathcal A}

\newcommand{\Da}{\mathscr{D}}
\newcommand{\Pa}{\mathscr{P}}
\newcommand{\one}{\mathbf{1}}

\newcommand{\lma}{\left\{}    
\newcommand{\rma}{\right\}}   

\newcommand{\lb}{\left(}
\newcommand{\rb}{\right)}
\newcommand{\lmb}{\left[}
\newcommand{\rmb}{\right]}
\newcommand{\lw}{\left|}
\newcommand{\rw}{\right|}

\newcommand{\ep}{\varepsilon}
\newcommand{\La}{\mathcal{L}}

\newcommand{\avg}[1]{\left\langle#1\right\rangle}

\newcommand\restr[2]{{%
  \left.\kern-\nulldelimiterspace
  #1
  \vphantom{\big|}
  \right|_{#2}
}}

\makeatletter
\renewcommand{\@@and}{ \ \ \ \ }
\makeatother

\title[Uniform boundary observability]
{Uniform boundary observability for periodic wave equations with low frequency exponent $4/3$}

\author{Zhonggan Huang}
\address{Institute for Theoretical Sciences, Westlake University,
  Hangzhou 310024, China}
\email{huangzhonggan@westlake.edu.cn}
\author{Jinping Zhuge}
\address{Morningside Center of Mathematics, Academy of Mathematics
  and Systems Science, Chinese Academy of Sciences, Beijing 100190, China}
\email{jpzhuge@amss.ac.cn}

\begin{document}

\begin{abstract}
We study boundary observability for Dirichlet wave equations with symmetric, periodic coefficients in bounded $C^{2,1}$ domains of dimension $d\ge2$. We establish two-sided observability estimates, uniform in the oscillation scale $\varepsilon$, for initial data generated by eigenfunctions with eigenvalues at most $\delta\varepsilon^{-4/3}$ with some small factor $\delta $. The proof reduces wave observability to boundary gradient estimates for narrow spectral clusters. These estimates are obtained using a second-order expansion and a uniform commutator estimate for harmonic extension in boundary layers.

\end{abstract}

\maketitle

\section{Introduction}

In this paper we study uniform boundary observability for the wave operator
\begin{equation}\label{intro:operator}
 \partial_t^2+\Lop_\eps,
 \qquad \Lop_\eps=-\operatorname{div}\big(A(x/\eps)\nabla\big),
\end{equation}
in a bounded domain $\Omega\subset\R^d$, $d\ge2$.
The parameter $\eps>0$ represents the scale of the periodic structure. We
assume that the matrix $A=(a_{ij})_{d\times d}$ is real, symmetric and
$\Z^d$-periodic, and satisfies
\begin{equation}\label{intro:ellipticity}
 \Lambda^{-1}|\xi|^2\le A(y)\xi\cdot\xi\le\Lambda|\xi|^2
 \qquad\text{for all }y,\xi\in\R^d,
\end{equation}
where $\Lambda\ge1$. In addition, we assume that $A$ is Lipschitz with the following gradient bound
\begin{equation}\label{intro:Lipschitz-A}
    |\nabla A(y)| \le \Lambda \qquad \text{for a.e. } y\in \R^d.
\end{equation}

We are interested in boundary observability, with constants independent of $\eps$, for initial data generated by the Dirichlet eigenfunctions of $\Lop_\eps$ with eigenvalues at most $\delta\eps^{-l}$, for some $l > 0$.  
The main challenge is to increase the exponent $l$ and thus enlarge the frequency bandwidth.

Let $\{(\lambda_{\eps,k},\psi_{\eps,k})\}_{k\ge1}$ be the Dirichlet
eigenpairs of $\Lop_\eps$ in $\Omega$, where the eigenvalues are listed in
nondecreasing order, with multiplicity, and the eigenfunctions form an
orthonormal basis of $L^2(\Omega;\mathbb C)$. For $N=N_\eps>0$, set
\begin{equation}\label{intro:spectral-space}
 \mathcal A_\eps^N
 =\operatorname{span}\{\psi_{\eps,k}:\lambda_{\eps,k}\le N\}.
\end{equation}
We consider the initial-Dirichlet problem
\begin{equation}\label{intro:wave} 
 \begin{cases}
  (\partial_t^2+\Lop_\eps)u_\eps=0&\text{in }\Omega\times(0,T),\\
  u_\eps=0&\text{on }\partial\Omega\times(0,T),\\
  u_\eps(\cdot,0)=f_\eps,\quad
  \partial_tu_\eps(\cdot,0)=g_\eps&\text{in }\Omega.
 \end{cases}
\end{equation}
Its conserved energy is given by
\begin{equation}\label{intro:energy}
 E_\eps(t)=\frac12\int_\Omega
 \big\{A(x/\eps)\nabla u_\eps\cdot\overline{\nabla u_\eps}
                  +|\partial_tu_\eps|^2\big\}\dd x.
\end{equation}
The question is how large $N=N_\eps$ may be chosen so that the boundary
observation of every solution with data in
$\mathcal A_\eps^N\times\mathcal A_\eps^N$ is uniformly comparable to its
energy.

\begin{theorem}\label{intro:main}
Let $\Omega\subset\R^d$ be a bounded $C^{2,1}$ domain, and let $A$ satisfy the
assumptions in \eqref{intro:ellipticity} and \eqref{intro:Lipschitz-A}. There exist positive constants $\delta$, $T_0$, $\eps_0$ and $C$
such that, if
\begin{equation}\label{intro:cutoff}
 0<\eps\le\eps_0,\qquad N_\eps\le\delta\eps^{-4/3},\qquad T\ge T_0,
\end{equation}
then every solution of \eqref{intro:wave} with
$(f_\eps,g_\eps)\in\mathcal A_\eps^N\times\mathcal A_\eps^N$ satisfies
\begin{equation}\label{intro:observation}
 C^{-1}E_\eps(0)
 \le \frac1T\int_0^T\int_{\partial\Omega}|\nabla u_\eps|^2\dd\sigma\dd t
 \le C E_\eps(0).
\end{equation}
The constants are independent of $\eps$, $N$, $T$, and the initial data;
they may depend on $d$, $\Omega$, and $\Lambda$ in \eqref{intro:ellipticity} and \eqref{intro:Lipschitz-A}.
\end{theorem}

In the pioneering work \cite{lin_uniform_2021, lin_erratum_2022}, Lin and Shen proved uniform boundary observability with the following spectral cutoff
\begin{equation}\label{intro:LS}
 N_\eps \le c_0T^{-1}\eps^{-1/2},\qquad T\ge T_0.
\end{equation}
Theorem~\ref{intro:main} enlarges the established eigenvalue range from $\eps^{-1/2}$ to $\eps^{-4/3}$, and removes the dependence of the factor on the length of time $T$. The optimal scale $\delta\eps^{-2}$ as conjectured in \cite{lin_uniform_2021} remains open. Our result does not reach this scale, and no optimality of the exponent $4/3$ is asserted.   We would also like to mention that when we were finishing our manuscript, the $\eps^{-1}$ cutoff was obtained independently in \cite{geng_uniform_2026}.

The main new ingredient of the proof is a boundary gradient estimate for narrow clusters of Dirichlet eigenfunctions up to spectral cutoff $\eps^{-4/3}$ to the operator $\La_\eps$ in $\Omega$.  To state the cluster estimate, we denote $P_\eps=\Lop_\eps^{1/2}$ and define for $s>0$ and $\rho>0$
\begin{equation}\label{intro:cluster-space}
\begin{split}
      X_{s,\rho}
   & =\mathcal A_\eps^N\cap\ \textup{Ran} \
    \mathbf1_{[s-\rho,s+\rho]}(P_\eps)\\
&= \operatorname{span}\left\{\psi_{\eps,k}:\lambda_{\eps,k}\le N_\eps\textup{ and }\sqrt{\lambda_{\eps,k}}\in [s-\rho,s+\rho]\right\}.
\end{split}
  \end{equation}
Note that here the operator $P_\eps$ has eigenvalues exactly $\sqrt{\lambda_{\eps,j}}$. Here formally the relation between the eigenvalue scale $\lambda$ and the frequency scale $s$ is $\lambda = s^2$.

\begin{theorem}\label{intro:cluster}
Under the assumptions of Theorem~\ref{intro:main}, there exist positive
constants $\delta$, $\rho_0$, $\eps_0$, $c$, and $C$ such that for
\[
0<\eps\le\eps_0, \ N_\eps\le\delta\eps^{-4/3} \ \textup{ and } \ 0<\rho\le\rho_0,
\]
the following estimate holds for all $z_\eps\in X_{s,\rho}$
\begin{equation}\label{intro:cluster-est}
 c s^2\Norm{z_\eps}{L^2(\Omega)}^2
 \le\int_{\partial\Omega}|\nabla z_\eps|^2\dd\sigma
 \le C s^2\Norm{z_\eps}{L^2(\Omega)}^2.
\end{equation}
The constants $C, c > 0 $ are independent of $\eps$, $N$, $s$, $\rho$, and
$z_\eps$.
\end{theorem}

The passage from Theorem~\ref{intro:cluster} to
Theorem~\ref{intro:main} uses a spectral localization argument in the
spirit of Miller's resolvent criteria \cite{miller_controllability_2005, Miller2012}. We prove
an abstract form in Section~\ref{sec.abstract.reduction} for future use. We also remark that the abstract framework we construct in Section~\ref{sec.abstract.reduction} is robust in all choices of frequencies, that is to say, the passage from Theorem~\ref{intro:cluster} to
Theorem~\ref{intro:main} is independent of the low frequency cut-off $\lambda\in [0,N_\eps)$, and in particular, the same framework holds also for the high frequency case $\lambda \in [N_\eps' , \infty)$. We refer to Remark \ref{r.generality.of.application} for more details.

Theorem~\ref{intro:cluster} improves the boundary gradient estimates for Dirichlet eigenfunctions in \cite{kenig_estimates_2013} from the range $\eps^{-1}$ to $\eps^{-4/3}$, which is of independent interest. However, we do not claim the optimality of the eigenvalue cutoff $\eps^{-4/3}$ either. The critical threshold is also $\eps^{-2}$. In fact, counterexamples in \cite{castro_low_2000} with $ d =1$ shows that eigenpairs can be constructed such that $\lambda_\eps \simeq \eps^{-2}$ and
\begin{equation*}
    \Norm{\nabla \psi_\eps}{L^2(\partial \Omega)} \simeq \lambda_\eps^{3/4}.
\end{equation*}
We refer to \cite{kenig_estimates_2013} for a survey on the asymptotic behavior of eigenvalue problems.

\subsection*{Techniques in the proofs} Let us briefly introduce some techniques in the paper.

\subsubsection*{Abstract reduction to narrow eigenclusters}

The abstract reduction framework in Section~\ref{sec.abstract.reduction} reduces the proof of Theorem \ref{intro:main} to that of Theorem \ref{intro:cluster} alone. Therefore, we refer to Section~\ref{sec.abstract.reduction} for a complete description of the reduction, and mainly introduce the techniques in proving the latter theorem on eigenclusters.

\subsubsection*{Second-order homogenization expansion}
One of the main ideas in the proof is the second order expansion of a unit vector $z_\ep \in X_{s,\rho}$, which is the main approach for building fine estimates for $z_\eps$ and the gradients near the boundary. We write $\La_0$ as the homogenized operator of $\La_\eps$ as $\eps\rightarrow 0$, and then solve for the unique $v\in H_0^1(\Om)$ of the following PDE
\[
\La_0 v = \La_\ep z_\eps = s^2 z_\eps + r_\eps, \textup{ with }r_\eps:= (\La_\ep - s^2) z_\eps.
\]
Here $r_\eps$ is small in the sense that $\norm{r_\eps}_2 \le C \rho s$ as $z_\eps$ lies in the narrow cluster $X_{s,\rho}$.

We can write $z_\eps$ as a second-order expansion as follows
\begin{equation}\label{intro.eq.expansion}
  z_\eps =  v+\eps \chi_j^\eps \pt_j v + \eps B_{1,\eps} + \eps^2 \mu_{ij}^\eps\pt_{ij}v + \eps^2 B_{2,\eps} + R_{2,\eps},
\end{equation}
where $\chi_j^\eps$ and $\mu_{ij}^\eps$ are respectively the first- and second-order correctors, and $B_{i,\eps}$ the $i$-th boundary layer that cancels the boundary values of $\eps \chi_j^\eps \pt_j v$ and $\eps^2 \mu_{ij}^\eps\pt_{ij}v$ respectively. We refer to \eqref{eq.second.order.ansatz} for a rigorous description of this expansion. The critical exponent can be heuristically seen from the following estimate for the remainder of the series above (we denote $\lambda=s^2$ as the scale of the clustered eigenfunctions)
\[
\norm{R_{2,\eps}}_{H^1(\Om)} \le C \eps^2 \norm{v}_{H^3(\Om)} \le C \eps^2 \lambda \norm{z_\eps}_{H^1(\Om)} \le C \eps^2 \lambda^{3/2},
\]
which becomes small when $\lambda \le N_\eps \le \delta \eps^{-4/3}$ for a small factor $\delta>0$, where $\lambda = s^2$ represents the scale of the eigenvalues of the clustered eigenfunctions.

We remark here that the second order expansion improves the $\eps^{-1}$ result in \cite{geng_uniform_2026}, which mainly uses the first order expansions, but it seems impossible to improve the exponent $4/3$ by simply increasing the order of expansion in \eqref{intro.eq.expansion}. We refer to Remark \ref{r.barrier.of.higher.order.expansion} for more details about this issue.

\subsubsection*{Boundary layer and commutator estimate} The entire path of proving \eqref{intro:cluster-est} is illustrated as follows
\begin{equation*}
    \Norm{\nabla v}{L^2(\partial \Omega)}\xRightarrow{A} \Norm{\nabla v}{L^2( \Omega_{\eps})} \xRightarrow{B} \Norm{\nabla z_\eps}{L^2( \Omega_{\eps})}\xRightarrow{C} \Norm{\nabla z_\eps}{L^2(\partial \Omega)}.
\end{equation*}
We first establish the boundary estimate for $v$ on the leftmost of the above display, in which we take the advantage that $\La_0$ is an operator with constant coefficients, and hence there is no extra blowing up terms in the Rellich identity. The passage $A$ is obtained through a modified trace theorem. The passage $B$ is one of the main ingredient, and it requires an estimate of the form
(see Lemma~\ref{l.grad.error.z.v}) 
\begin{equation}\label{intro:collar-error}
 \Norm{\nabla z_\eps - (\nabla \Phi_\eps) \nabla v}{L^2(\Omega_{b\eps})}
 \le C_b\big\{
 \eps^{3/2}\Norm{v}{H^2(\Omega)}^{1/2}
             \Norm{v}{H^3(\Omega)}^{1/2}
 +\eps^2\Norm{v}{H^3(\Omega)}\big\},
\end{equation}
where $\Phi_\eps$ is the Dirichlet corrector and $\nabla  \Phi_\eps$ is an invertible matrix in $\Omega_{b\eps}$ for small $b>0$ (see \cite[Lemma~5.1]{kenig_estimates_2013}). The passage $C$ follows from the classical good-bad patch argument and unique continuation principle in \cite[Theorem~1.4]{kenig_estimates_2013}.

We emphasize that the proof of \eqref{intro:collar-error} is a major difficult step. The difficulty mainly comes from the boundary layer term $B_{1,\eps}$, which involves a new commutator estimate for $\Lop_\eps$-harmonic extension.
Write $\mathcal H_\eps$ for the Poisson extension associated with
$\Lop_\eps$ (see \eqref{eq:Poisson-integral} for a precise definition), and put
\[
 [\mathcal H_\eps,\phi]b
 =\mathcal H_\eps(\phi b)
                  -\phi \mathcal H_\eps(b),
\]
where $b_j\in L^\infty(\partial \Omega)$ and $\phi_j \in H^1(\Omega)$. In our application $b = \chi_j(x/\eps)$ and $\phi = \partial_j v$, and thus $B_{1,\eps} = \mathcal{H}_\eps (\phi_j b_j)$.
Let $\Omega_t=\{x\in\Omega:\operatorname{dist}(x,\partial\Omega)<t\}$.
For every fixed $B>0$, we establish, for sufficiently small $\eps$,
\begin{equation}\label{intro:commutator}
    \begin{split}
    \Norm{\nabla[\mathcal H_\eps,\phi]b}{L^2(\Omega_{B\eps})}
    \le C_B\Norm{b}{L^\infty(\partial\Omega)}
    \big\{\eps^{1/2}\Norm{\phi}{L^2(\partial\Omega)}
       +\Norm{\nabla\phi}{L^2(\Omega_{2B\eps})}\big\}.
    \end{split}
\end{equation}
Here no derivative of the bounded boundary
function $b$ occurs, which is essential when $b$ is an oscillating
corrector in our application. This estimate is of independent interest and is proved in
Section~\ref{sec.sharp.commutator.estimate} for $C^{2,1}$ domains
by Poisson-kernel expansions
of Kenig, Lin, and Shen \cite{KLS14}, together with the
Dirichlet-to-Neumann commutator estimates in
\cite[Theorem~5.1]{KLS14} and \cite[Theorem~1.2]{Shen16}.

\subsection*{Literature}
Boundary observability for fixed coefficients is classical: the Hilbert uniqueness method
of Lions \cite{LionsSIAMRev1988,Lions} and the geometric control condition of Bardos,
Lebeau, and Rauch \cite{BardosLebeauRauch1992}, quantified by Laurent and L\'eautaud
\cite{laurent_uniform_2016}. For oscillating coefficients, Lions
\cite{lions_controlabilite_1988} related optimal controls to those of the homogenized
problem, and Avellaneda and Lin \cite{AvellanedaLinPoisson}, using their uniform elliptic
estimates \cite{AvellanedaLin}, described the weak limits of boundary controls.
Avellaneda, Bardos, and Rauch \cite{Avellaneda1992ContrllabilitEH} then showed that
$\varepsilon$-uniform estimates fail on the full spectrum, even in one dimension. This motivates
the restriction to $\mathcal{A}_\eps^N$ for the initial data.

For the uniform control in homogenization setting,
in one dimension, Castro \cite{castro_boundary_1999} proved sharp partial controllability
on the band $\lambda\le\delta\varepsilon^{-2}$ and on a high-frequency band, with spectral
asymptotics from Castro and Zuazua \cite{castro_low_2000,castro_high_2000} and related
results of Pedregal and Periago \cite{pedregal_remarks_2006}. In higher dimensions, Lebeau
\cite{lebeau_wave_2000} reached $\varepsilon^{-2}$ for interior observation and
density-type coefficients. For divergence-form boundary observation, Lin and Shen
\cite{lin_uniform_2021,lin_erratum_2022} obtained $N\le c_0T^{-1}\varepsilon^{-1/2}$. A recent preprint of 
Geng-Wang-Wang \cite{geng_uniform_2026} obtained $N\le\delta\varepsilon^{-1}$.
Our Theorem~\ref{intro:main} gives $N\le\delta\varepsilon^{-4/3}$; the endpoint $\varepsilon^{-2}$ remains
open.
 
Boundary estimates for eigenfunctions and clusters
are due to Hassell and Tao \cite{hassell_upper_nodate} (see the errata
\cite{hassell_errata_nodate}), Barnett and Hassell \cite{BarnettHassell2011}, and Xu
\cite{xu_upper_2012} in the smooth-coefficient setting.   For periodically oscillating coefficients, Kenig-Lin-Shen
\cite{kenig_estimates_2013} obtained the estimates for eigenvalues having scale $\lambda\le C\varepsilon^{-1}$. Theorem~\ref{intro:cluster} extends the two-sided cluster bound to
$\varepsilon^{-4/3}$, using a second-order adaptation of Kenig-Lin-Shen
\cite{KLS14} and the commutator bounds of Shen \cite{Shen16}; the abstract wave reduction follows Miller
\cite{miller_controllability_2005,Miller2012}.

\subsection*{Notation}
{Throughout this paper, all the constants $C, c, B, b$ etc. are independent of $\eps$ and $\lambda$. They may vary from line to line. The Einstein summation convention is applied for repeated indices.}

\subsection*{Outline of the paper}
The paper is organized as follows. Section~\ref{sec.abstract.reduction} establishes the abstract
spectral reduction. Section~\ref{sec.body.estimate.43} develops the second-order expansion and
the boundary estimates used in Theorems~\ref{intro:main} and
\ref{intro:cluster}. Section~\ref{sec.sharp.commutator.estimate} proves the sharp commutator
estimate. Auxiliary regularity and boundary-layer estimates are collected
in the appendix.

\subsection*{AI use declaration}
The mathematical results in this paper are obtained through multiple rounds of discussion with OpenAI ChatGPT 5.6 Sol. The paper itself is, however, completely written by the human authors after careful review and reorganizations of the AI outcomes. Both authors are fully responsible for the accuracy and integrity of the manuscript.

\subsection*{Acknowledgments}
J.~Z. is partially supported by NNSF of China (Nos. 12494541, 12288201, 12471115). Z.~H. acknowledges the helpful discussions with Zhongwei Shen.

\section{Reduction to clustered eigenfunctions}
\label{sec.abstract.reduction}

In this section, we discuss a reduction from wave equation controllability to the boundary gradient estimates of clustered eigenfunctions. This reduction adapts from the work of Miller \cite{miller_controllability_2005}. A similar form of this framework also appeared in a recent work of Geng-Wang-Wang \cite{geng_uniform_2026}.

\subsection{Spectral decomposition of self-adjoint operators}

In this subsection we recall some classical results on the spectral properties of an abstract self-adjoint operator. On a Hilbert space $W$, a self-adjoint operator $J$ satisfies
\[
(Jx,y)=(x,Jy) \textup{ for all }x,y\in \Da(J) \subset W,
\]
where $ \Da(J)$ is a dense subset of $W$ on which $J$ is defined.

On the other hand, the spectral theorem decomposes $J$ into a projection-valued measure. To be more precise, a projection-valued measure $\mathscr{P}$ is a Borel measure on $\R$ that takes values in the projections of $W$, that is for every Borel $I\subset \R$, $\Pa(I)$ is a projection on $W$.

\begin{theorem}[Spectral theorem]
    \label{t.spectral.theorem}
    To every self-adjoint operator $J$ there is a unique projection-valued measure $\Pa_J$ such that 
    \begin{equation}
        \label{eq.spectral.theorem}
        J = \int_\R \lambda \dd \Pa_J(\lambda).
    \end{equation}
\end{theorem}

This theorem is a classical result in the spectral theory, see  \cite[Theorem 3.6]{teschl_mathematical_2014}. With this theorem one can easily define for any Borel function $f: \R \rightarrow \mathbb C$, the following operator
\begin{equation}
    \label{eq.define.fs}
    f(J) : = \int_\R f(\lambda) \dd \Pa_J(\lambda).
\end{equation}
Here, one can show that the domain of $f(J)$ is dense in $W$. Common choices of $f$ are $\one_I(J)$ for some Borel set $I\subset \R$ and $e^{it J}$.

\subsection{Spectral  cut-off}

Let $L$ be a self-adjoint operator on a complex Hilbert space $\hat{X}$, with domains denoted as $\Da(L)\subset \hat X$. We assume that $L$ is positive in the sense that there is a constant $p_0>0$ such that for all $z\in \Da(L)$
\begin{equation}
    \label{eq.coercive.P.alpha}
    (L z, z) \ge p_0^2 \|z\|_{\hat X}^2. 
\end{equation}

We denote for any Borel subset $I \subset \R$
\begin{equation}
    \label{eq.half.operator}
    P : = L^{1/2} \textup{ and }X : = \textup{Ran} ~ \one_{I} (L).
\end{equation}

\begin{remark}
    \label{r.generality.of.application}
    Here $I$ can be chosen as $[0,N)$ for low frequency and $(N',\infty)$ for high frequency. We remark here that the abstract framework we establish in this section is independent of the choice of low or high frequency. The only ingredients that matter are the bounding constants that will show up in \eqref{eq.cluster.estimate} and \eqref{eq.uniform.control} below. 
\end{remark}

Let $B: X\rightarrow Y$ be a linear operator with $Y$ an arbitrary Hilbert space. We require 
\begin{equation}
    \label{eq.continuity.Ba}
    \norm{B z}_{Y} \le b \lb \norm{z}_{X} + \norm{L z}_{X}\rb, \ \textup{ and }\ \Da(B) = X \cap \Da(L) =:X_L,
\end{equation}
where $0<b<\infty$ is an arbitrary number.

The last assumption we need is the following quantitative packet condition: for all $s\in \R$
\begin{equation}
    \label{eq.cluster.estimate}\tag{S}
 m \|P z\|_{X}^2 \le \|B z\|_{Y}^2 \le M\|P z\|_{X}^2, \qquad \textup{ for }z\in X \cap  \textup{Ran} ~ \one_{[s-\rho,s+\rho]} (P)
  \end{equation}
where $\rho,m,M>0$ are constants independent of $s$. Note that it is not difficult to show that 
\[
X \cap  \textup{Ran} ~ \one_{[s-\rho,s+\rho]} (P) \subset  \Da(L)\subset\Da(P).
\]

\subsection{Abstract wave equation and uniform control}

In this subsection we follow the same notation as in the above subsection. We discuss the solution $u(t)$ to the following wave equation
\begin{equation}
    \label{eq.wave.abstract}
    \begin{cases}
        u_{tt} + L u = 0 & \textup{ for }t>0\\
        u(0) = f \textup{ and }u_t(0)=g & \textup{ where }f\in X_P:= X\cap \Da(P) \textup{ and }g\in X.
    \end{cases}
\end{equation}

Note that the solution $u$ can be written explicitly as 
\begin{equation}
    \label{eq.wave.abstract.sol}
    u(t) = \frac{e^{it P}+ e^{-it P}}{2} f + \frac{e^{it P}- e^{-it P}}{2i} P^{-1} g.
\end{equation}
Observe that because $(f,g)\in X_P\times X$ (here $X_P : = X \cap \Da(P)$), we have for each $t$
\begin{equation}
    \label{eq.wave.abstract.sol.regularity}
    u(t) \in X_P , \ P  u(t) \in X \textup{ and } u_t (t) \in X.
\end{equation}
Note that the vector $Bu(t)$ at this stage is well-defined only when
\begin{equation}
    \label{eq.well-define.But}
    (f,g) \in X_L\times X_P \textup{ and then }u(t) \in X_L \textup{ for each }t.
\end{equation}
That is to say, $B$ can be regarded as a linear operator $B^*$ from $X_L\times X_P$ to $L^2([0,T];Y)$ with the form
\begin{equation}
    \label{eq.define.Bstar}
    B^*(f,g): = Bu(t) \in L^2([0,T];Y).
\end{equation}
The uniform controllability inequality for the solution $u$ with respect to the operator $B$ takes the following form: there are constants $C, T_0>0$ such that for $T\ge T_0$, the operator $B^*$ defined in \eqref{eq.define.Bstar} has a unique extension to $X_P\times X$, and 
\begin{equation}
    \label{eq.uniform.control}\tag{U}
    C^{-1} E\le \frac{1}{T}\norm{B^*(f,g)}_{L^2([0,T];Y)}^2 \le C E, \qquad E:= \|P f\|_{X}^2 + \|g\|_{X}^2.
\end{equation}

\begin{theorem}
    \label{t.SimpliesU}
  Assume \eqref{eq.coercive.P.alpha}, \eqref{eq.half.operator} and \eqref{eq.continuity.Ba}.  Let $u(t)$ be a solution to the wave equation \eqref{eq.wave.abstract} with initial data $$(u(0),u_t(0))=(f,g)\in X_P \times X,$$ then \eqref{eq.cluster.estimate} implies \eqref{eq.uniform.control}. Moreover, the constant $C$ in \eqref{eq.uniform.control} depends only on $p_0, m,M$ in \eqref{eq.cluster.estimate}, and the constant $T_0$ depends only on $m,M$ and $\rho$ in \eqref{eq.cluster.estimate}.
\end{theorem}

To prove this theorem we need the following auxiliary lemma.

\begin{lemma}
    \label{l.band.lemma}
    Let $D$ be a self-adjoint operator (not necessarily compact resolvent) on a Hilbert space $V$, and $H: \Da{(D)} \rightarrow W$ is another linear operator satisfying
    \begin{equation}
        \label{eq.continuity.C}
        \norm{H v }_{W} \le h \lb \norm{v}_{V} + \norm{D v}_V\rb, \ \textup{ for some }h<\infty \textup{ and all }v\in \Da(D).
    \end{equation}
    Suppose there are positive constants $m,M$ and $r$ such that
    \begin{equation}
        \label{eq.conditionS.lemma.band}
        m\|z\|_V^2 \le \|H z\|_W^2 \le M \|z\|_V^2  \qquad \textup{ for }z\in\textup{Ran}~ \one_{[s-r,s+r]}(D)  \textup{ and all } s\in \R.
    \end{equation}
Then the operator $H e^{itD} : \Da(D) \rightarrow L^2([0,T];W)$ has a unique extension $H^*$ to whole $V$, and there is a constant $R>0$ and $T_0 = T_0(p_0, M,m,r) > 0$ such that for every $T\ge T_0$ and $x\in V$
\begin{equation}
    \label{eq.uniform.control.band}
  R^{-1} \|x\|_V^2 \le \frac{1}{T}\norm{H^*(x)}_{L^2([0,T];W)}^2 \le R \|x\|_V^2.
\end{equation}
Here the constant $R$ depends only on $m,M$ and $T_0$ depends only on $m,M$ and $r$. In particular, $R$ and $T_0$ are independent of the choice of the Hilbert spaces $X$ and $W$, the operator $D$ and the bounding constant $h$ in \eqref{eq.continuity.C}.

\end{lemma}

\begin{proof}

Let us first restrict our discussion to the space of bounded spectrum
\begin{equation}
    \label{eq.bounded.spectrum.space}
    V^K : = \textup{Ran} \ \one_{[-K,K]} (D)
\end{equation}
for some $K<\infty$. Note that $V^K\subset V \cap \Da(D) = \Da(H)$. Moreover, on $V^K$ both $D$ and $H$ are bounded operators even though the bounding constants may depend on $K$ and $h$ in \eqref{eq.continuity.C}.

Set $f(t) = \frac{\sin{(rt)}}{rt}$ with $f(0)=1$, where $r$ is the fixed band width in \eqref{eq.conditionS.lemma.band}. It is not difficult to derive that the Fourier transform $\hat{f}$ of $f$ has the form
\begin{equation}
    \label{eq.fouriertransformintime}
    \hat{f}(s) = \frac{1}{\sqrt{2\pi}}\int_\R f(t) e^{-its}\dd t = \sqrt{\frac{\pi}{2r^2}} \mathbf{1}_{[-r,r]}(s), \qquad \textup{and }\|f\|_{L^2(\R)}^2 = \pi/r.
\end{equation}

By applying Theorem \ref{t.spectral.theorem} and Fubini theorem, we obtain for $x\in V^K$ and Lebesgue almost every $s\in \R$
\begin{equation}
    \label{eq.fourier.wave.time}
    \begin{split}
       \widehat{f(t) e^{itD}} x(s) &= \frac{1}{\sqrt{2\pi}} \int_{\R} \lmb \int_{\R} e^{-it(s-\lambda)}f(t) \dd t\rmb \dd \Pa_D(\lambda)x  \\
       &=\hat{f}(s-D) x \in \textup{Ran}~ \one_{[s-r,s+r]}(D).
    \end{split}
   \end{equation}
On one hand, by applying the Plancherel identity on $L^2(\R;W)$ there is 
\begin{equation}
    \label{eq.plancherel}
    \int_\R \norm{H \hat{f}(s-D) x}_W^2 \dd s = \int_\R |f(t)|^2 \norm{H e^{itD} x}_W^2 \dd t
\end{equation}
On the other hand, by similarly applying the Plancherel identity on $L^2(\R;V)$
\begin{equation}
    \label{eq.norm.compute}
    \int_\R \norm{\hat{f}(s-D) x}_V^2 \dd s = \norm{f}_{L^2(\R)}^2 \norm{x}_V^2.
\end{equation}
By combining \eqref{eq.conditionS.lemma.band}, \eqref{eq.fouriertransformintime}, \eqref{eq.plancherel} and \eqref{eq.norm.compute}, we obtain
\begin{equation}
    \label{eq.fweight.uniform.control.band}
    \frac{m\pi}{r} \|x\|_V^2 \le \int_{\R} |f(t)|^2 \|H e^{it D} x\|_W^2 \dd t \le \frac{M\pi}{r} \|x\|_V^2 \ \textup{ for all }x\in V^K
\end{equation}
Note that the Fourier transform commutes with $H$ on $V^K$, where $H$ is fully defined and bounded. The main advantage of \eqref{eq.fweight.uniform.control.band} is that the possibly large norm of $H$ does not enter the estimate.
We also remark here that the same estimate in \eqref{eq.fweight.uniform.control.band} also holds for $f(t)$ being replaced by $f(t-t_0)$ for any fixed $t_0\in \R$. 

Now, we observe that because $\frac{\sin{(u)}}{u} \ge \sin(1)$ for all $|u|\le 1$, there is by \eqref{eq.fweight.uniform.control.band}
\begin{equation}
    \label{eq.upperbound.unit.interval}
    \int_{t_0-1/r}^{t_0+1/r} \|H e^{it D} x\|_W^2 \dd t \le  C\int_{\R} |f(t-t_0)|^2 \|H e^{it D} x\|_W^2 \dd t \le \frac{CM\pi}{r} \|x\|_V^2 
\end{equation}
for all $t_0 \in \R$ and $x\in V^K$. Integrating over $t_0 \in [0,T]$ gives the upper bound in \eqref{eq.uniform.control.band}.

In the following we use a constant $R$ that depends at most on $m,M$ but may be different from line to line. To prove the lower bound, we observe that the function $$|f(t)|\le\min\lma 1, \frac{1}{r |t|}\rma,$$ and therefore by \eqref{eq.upperbound.unit.interval}
\begin{equation}
    \label{eq.error.ext.interval}
    \begin{split}
         \int_{|t| \ge R/r} |f(t)|^2 \|H e^{itD} x\|_W^2 \dd t &\le \sum_{|j|\ge R-1} \frac{1}{j^2}\int_{j/r}^{(j+2)/r} \|H e^{it D} x\|_W^2 \dd t\\[5pt]
         &\le \frac{R\|x\|_V^2}{r} \sum_{|j|\ge R-1} j^{-2} ,
    \end{split}
   \end{equation}
where the last inequality is obtained by combining \eqref{eq.upperbound.unit.interval}. By choosing $R$ sufficiently large that is only depending on $m,M$, we obtain 
\begin{equation}
    \label{eq.error.ext.interval.done}
    \int_{|t| \ge R/r} |f(t)|^2 \|H e^{itD} x\|_W^2 \dd t \le \frac{m\pi}{2r} \|x\|_V^2.
\end{equation}
This shows that by combining \eqref{eq.fweight.uniform.control.band}
\begin{equation}
    \label{eq.lowerbound.unit.interval}
    \int_{-R/r}^{R/r} \|H e^{itD} x\|_W^2 \dd t \ge \frac{m\pi}{2r} \|x\|_V^2.
\end{equation}
Thus the estimate \eqref{eq.uniform.control.band} is obtained for $x\in V^K$ by observing that the time interval $[0,T]$ is covered by $O(rT)$ intervals of length $2/r$, and $O(rT/R)$ intervals of length $2R/r$. Also note that the proof requires the minimal time of observation $T_0 = 2R/r$.

To finish the proof one needs to extend the definition of $H e^{it D}$ to whole $V$. Currently it is only defined on $\cup_{K>0} V^K$, which is a dense subset of $V$ because 
\[
V^K \ni  \one_{[-K,K]}(D) x \rightarrow x
\]
for all $x\in V$ as $K \rightarrow \infty$. Because the bounding constants $R$ and $T_0$ in \eqref{eq.uniform.control.band} are independent of $K$, the operator $H e^{it D}$ can be  uniquely extended to be a bounded linear operator $\Tilde{H}_t : V \rightarrow L^2([0,T]; W)$. Here $\Tilde{H}_t$ coincides with $H e^{itD}$ on $\Da(D)$ because of the assumption \eqref{eq.continuity.C} and the unitarity of $e^{itD}$.

\end{proof}

\begin{proof}[Proof of Theorem \ref{t.SimpliesU}]
    The solution $u(t)$ to the wave equation \eqref{eq.wave.abstract} has the following explicit form
\begin{equation}
    \label{eq.explicit.wave.sol}
    \begin{split}
         u(t) &= \frac{e^{it P} + e^{-it P}}{2} f + \frac{e^{itP}-e^{-it P}}{2i} P^{-1} g\\
         &=P^{-1}\lb e^{it P} x_+ + e^{-it P} x_- \rb,
    \end{split}
\end{equation}
where 
\begin{equation}
    \label{eq.define.xpm}
    x_+ = \frac{P f - i g}{2} \textup{ and }x_- = \frac{P f + i g}{2}.
\end{equation}

This shows that for $(f,g)\in X_L \times X_P$, the generated vector $(x_+,x_-)\in X_P\times X_P$. Now we denote $V : = X \times X$ and $W:= Y$. We define on $V$ the following operators
\begin{equation}
    \label{eq.apply.abstract.D}
    D : = \lb\begin{matrix}
        P & 0\\
        0 & -P
    \end{matrix}\rb \textup{ and }H(x_1,x_2) := BP^{-1}(x_1+x_2).
\end{equation}
We will apply Lemma \ref{l.band.lemma} to $(D,H)$. First we note that $P^{-1}: X_P \to X_L$ is well-defined and by the assumption \eqref{eq.continuity.Ba}, for $z = (z_1,z_2) \in V$,
\begin{equation}
\begin{aligned}
    \norm{Hz}_{W} = \norm{BP^{-1}(z_1+z_2)}_{W} & \le b \lb \norm{P^{-1}(z_1 + z_2)}_X + \norm{P(z_1 + z_2)}_X) \rb \\
    & \le  \sqrt{2} b\lb p_0^{-1} \norm{z}_V + \norm{Dz}_V \rb,
\end{aligned}
\end{equation}
where we have used $\norm{P^{-1}z_1}_X \le p_0^{-1} \norm{z_1}_X$, due to \eqref{eq.coercive.P.alpha}.
This verifies the condition \eqref{eq.continuity.C} in Lemma \ref{l.band.lemma} with $h = \sqrt{2}b \max\{p_0^{-1}, 1\}$.

By requiring the band width $r = \min \{\rho, p_0 /2\}$, we obtain that an interval of length $2r$ cannot intersect the spectrum of both $P$ and $-P$. Thus for any $s\in \R$, if $z = (z_1, z_2) \in \textup{Ran} ~\one_{[s-r,s+r]}(D)$, then either $z_1 = 0$ or $z_2 = 0$. Suppose without loss of generality, $z_2 = 0$. Then 
\[
Hz = BP^{-1}z_1, \ \norm{z}_V = \norm{z_1}_2  \quad \textup{and} \quad \|Hz\|_{W}^2 = \|BP^{-1}z_1\|_{W}^2.
\]
By the assumption \eqref{eq.cluster.estimate} 
\begin{equation}
    \label{eq.apply.abstract}
    m \|z\|_V^2 \le \|Hz\|_{W}^2 \le M \|z\|_V^2
\end{equation}
for all $z\in V \cap \textup{Ran}~ \one_{[s-r,s+r]}(D)$ and $s \in \R$. This verifies \eqref{eq.conditionS.lemma.band}.

Hence, Lemma \ref{l.band.lemma} yields for any $z\in V$
\begin{equation}\label{eq.Abstract-observe}
    R^{-1} \|z\|_V^2 \le \frac{1}{T}\norm{H e^{itD} }_{L^2([0,T];W)}^2 \le R \|z\|_V^2.
\end{equation}
Let $z = (x_+, x_-) \in V$.
By \eqref{eq.explicit.wave.sol} and \eqref{eq.apply.abstract.D}, we have $Bu(t) = He^{itD}(x_+, x_-)$. Hence, \eqref{eq.uniform.control} follows from \eqref{eq.Abstract-observe} and 
\[
E = \|P f\|_{X}^2 + \|g\|_{X}^2  =2(\|x_+\|_X^2+\|x_-\|_X^2) = \norm{z}_V^2,
\]
where the last identity follows from 
\eqref{eq.define.xpm}.
\end{proof}

\section{A uniform bound for spectral cutoff \texorpdfstring{$ \eps^{-4/3}$}{sdas}}
\label{sec.body.estimate.43}

We establish the proof of the uniform bound \eqref{intro:observation} and \eqref{intro:cluster-est} with $N_\eps \simeq \eps^{-4/3}$. We use the formalism established in Section \ref{sec.abstract.reduction}. Specifically, we use the following notations:
\begin{itemize}
    \item We denote $$L = \Lop_\eps = - \partial_i (a_{ij}(x/\eps)\partial_j)$$ as the self-adjoint linear operator on the Hilbert space $\hat X := L^2(\Omega;\mathbb C)$ with the domain of definition $$\Da(L) = H^2\cap H_0^1(\Om).$$  
    For convenience we denote $\norm{\cdot}_2 : = \norm{\cdot}_{\hat X}$.
    \item Note that even though we work in complexified spaces in Section \ref{sec.abstract.reduction}, the original controllability problem is for real-valued functions. Therefore for simplicity we only need to establish estimates for real-valued functions in this section and do not write down the conjugate sign for inner product, i.e.
    \[
    (u,v) = \int_{\Om} u \Bar{v} \dd x = \int_{\Om} uv \dd x
    \]
    for all $u ,v $ being real functions. 
    \item The operator $\La_\eps$ has a compact resolvent $\Lop_\e^{-1}$ on $L^2(\Omega)$, which is bounded into $H^1_0(\Omega)$. Hence it has a sequence of eigenpairs $(\lambda_{\eps,j}, \psi_{\eps,j})$ such that $\lambda_{\eps,j}$ is nondecreasing in $j$ and for all $j\ge 1$ 
    $$\La_\eps \psi_{\ep,j} = \lambda_{\eps,j}\psi_{\ep,j}.$$
   The eigenfunctions $\{\psi_{\eps,j}\}$ form an orthonormal basis for $\hat X = L^2(\Omega;\mathbb C)$.
  \item Denote $P=P_\eps : = \Lop_\eps^{1/2}$. In $\Da(\Lop_\eps)$, it holds $\Lop_\eps = P_\eps^2$. Moreover, for $f,g \in H^1_0(\Omega)$,
  \begin{equation}\label{eq.Pe-integral-form}
  (P_\eps f, P_\eps g) = \int_{\Omega} A(x/\eps) \nabla f \cdot \overline{\nabla g}. 
  \end{equation}
   Note that
    \begin{equation}
        \Da(P_\eps) = \left\{f = \sum_{j\ge 1} a_j \psi_{\eps,j}: \sum_{j\ge 1} \lambda_{\eps,j}|a_j|^2 < \infty  \right\} = H_0^1(\Omega).
    \end{equation} 
  By the uniform ellipticity of $A = (a_{ij})_{d\times d}$ as in \eqref{intro:ellipticity}, we obtain by Poincar\'{e} inequality
    \begin{equation}
        \label{eq.coercive.43}
        \norm{P_\ep z }_2 \ge p_0\norm {z}_2 \textup{ for all }z\in H_0^1(\Om)\subset \hat{X}
    \end{equation}
    for some positive constant $p_0$ independent of $\eps$. This guarantees the coercivity condition required in \eqref{eq.coercive.P.alpha}.

    \item  We choose $I_\eps = [0,N_\eps]$, where $N_\eps = \delta \eps^{-4/3}$ with the small $\delta>0$ to be chosen. For this choice of $N_\eps$ we define
    \[
    X := \textup{Ran} \ \one_{I_\eps} (L) \subset \Da(L).
    \]
    Note that because $L = \Lop_\eps$ has compact resolvent, the space $X$ is finite dimensional. Also the space $X$ coincides with $\AN_\ep^N$ as defined in \eqref{intro:spectral-space}.
      
    \item Define $Y : = L^2(\pt \Omega; \mathbb C^d)$, and denote 
    \[
    \begin{split}
        B : X &\rightarrow Y\\
        u & \mapsto \restr{\gd u}{\pt \Omega}
    \end{split}
    \]
    for every $u \in X \cap \Da(L)=X$ (the last identity holds because we are working in the lower frequency regime).
    \item We define for $\lambda >0$ and $\rho>0$ the following spaces
        \begin{equation}
        \label{eq.define.X.s.rho.43}
        \begin{aligned}
            X_{\sqrt{\lambda},\rho} &:=\textup{Ran} \ \one_{[0,N_\eps)} (\Lop_\eps) \  \cap \  \textup{Ran} \ \one_{[\sqrt{\lambda}-\rho,\sqrt{\lambda}+\rho]} (P_\eps) \\
            & = \spn\left\{ \psi_{\eps,k}: \lambda_{\eps,k} \le N_{\eps}, \sqrt{\lambda} - \rho \le \sqrt{\lambda_{\eps,k}} \le \sqrt{\lambda} + \rho \right\}.
        \end{aligned}
    \end{equation}
    Note that $X_{\sqrt{\lambda},\rho} \subset \hat{X} = L^2(\Om;\mathbb C)$.
 
\end{itemize}

\subsection{Some basic estimates}

By the definition \eqref{eq.define.X.s.rho.43} of $X_{\sqrt{\lambda},\rho}$,  it is possible that $X_{\sqrt{\lambda},\rho}$ is empty. In the following, we automatically assume $X_{\sqrt{\lambda},\rho}$ has at least one eigenfunction. Then for $z_\eps \in X_{\sqrt{\lambda},\rho}$, we can write
\begin{equation}
    z_\eps = \sum a_k \psi_{\eps, k},
\end{equation}
where the sum is taken over all $k$ with $\lambda_{\eps,k} \le N_\eps$ and $|\sqrt{\lambda_{\eps,k}} - \sqrt{\lambda}| \le \rho$. It follows that
\begin{equation}\label{eq.LePe-formula}
    \Lop_\eps z_\eps = \sum a_k \lambda_{\eps,k} \psi_{\eps, k} \in X_{\sqrt{\lambda},\rho}, \qquad P_\eps z_\eps = \sum a_k \sqrt{\lambda_{\eps,k}} \psi_{\eps, k} \in X_{\sqrt{\lambda},\rho}.
\end{equation}

\begin{lemma}
    \label{l.estimates.Zeps.43}
  Assume that $\rho < \min\{p_0/3, 1/2\} $.  Let $z_\eps$ be an arbitrary unit vector in $X_{\sqrt{\lambda},\rho}$, then for all integer $k>0$
    \begin{equation}
        \label{eq.estimates.Peps.Zeps.43}
        2^{-k} \lambda^{k/2} \le \norm{P_\eps^k z_\eps}_{2} \le  2^k \lambda^{k/2} .
    \end{equation}

\end{lemma}

\begin{proof}
    Note that by \eqref{eq.LePe-formula}
    \[
  (\sqrt{\lambda}-\rho)^k \norm{z_\eps}_{2}\le  \norm{P_\eps^k z_\eps}_{2} \le  (\sqrt{\lambda}+\rho)^k \norm{z_\eps}_{2}.
    \]
    When $\rho < \min\{p_0/3, 1/2\} $, we obtain that either $\textup{Ran} \ \one_{[\sqrt{\lambda}-\rho,\sqrt{\lambda}+\rho]}(P_\eps)$ is trivial or $\sqrt{\lambda} \ge 2\rho$. The inequality \eqref{eq.estimates.Peps.Zeps.43} is obtained by plugging $\sqrt{\lambda} \ge 2\rho$ into the above display.
\end{proof}

For $z_\ep \in X_{\sqrt{\lambda},\rho}$, we are interested in the following auxiliary function $$v\in H_0^1(\Omega) = \Da(P_\eps)$$ satisfying
\begin{equation}
    \label{eq.homog.v.43}
        \Lop_0 v = \Lop_\eps z_\eps  = \lambda z_\eps + r_\eps, \quad \textup{with } r_\eps : = \Lop_\eps z_\eps  - \lambda z_\eps,
\end{equation}
where the operator $\Lop_0 = - \partial_i(\hat a_{ij} \partial_j)$ is the homogenization limit operator of $\Lop_\eps$ as $\eps \rightarrow 0$ by standard homogenization theory; see \eqref{eq.homogenized.tensor} below for a definition of the homogenized matrix $\widehat{A} = (\widehat{A}_{ij})$. Note that 
\begin{equation}
    \label{eq.bound.reps.43}
    \norm{r_\eps}_{2} =  \norm{(P_\eps - \sqrt{\lambda})(P_\eps+\sqrt{\lambda})z_\eps}_{2} \le \rho(2\sqrt{\lambda}+\rho)\norm{z_\eps}_{2}.
\end{equation}

\begin{lemma}
    \label{l.estimate.homog.v.43}
 Suppose $\Omega$ is $C^{2,1}$, $z_\eps\in X_{\sqrt{\lambda},\rho}$ is a unit vector and $\rho < \rho_0 $ for some small $\rho_0$ independent of $\eps$ and $\lambda$.   The auxiliary function $v\in H_0^1(\Omega)$ solving \eqref{eq.homog.v.43} satisfies the following estimates for a  constant  $C >0$
 \begin{equation}
     \label{eq.H2.homog.v.43}
     \norm{\gd^j v}_{2} \le  C\lambda^{j/2}, \quad j = 1,2,3.
 \end{equation}

\end{lemma}

\begin{proof}
    By \eqref{eq.Pe-integral-form} and \eqref{eq.estimates.Peps.Zeps.43}, $\norm{\nabla z_\eps}_2 \simeq \norm{P_\eps z_\eps} \simeq \lambda^{1/2}$. Testing the equation $\Lop_0 v = \Lop_\eps z_\eps$ by $v$, we can obtain 
    \begin{equation}
        \norm{\nabla v}_2^2 \le C (\La_0 v, v)= C(\La_\ep z_\ep, v) \le C \norm{\nabla z_\eps}_2 \norm{\nabla v}_2,
    \end{equation}
    yielding \eqref{eq.H2.homog.v.43} with $j = 1$. By \eqref{eq.LePe-formula}, we see $\norm{\Lop_\eps z_\eps}_2 \le C \lambda$.
    The $H^2$ regularity estimate applied to $\Lop_0 v = \Lop_\eps z_\eps$ in a bounded $C^{1,1}$ domain implies \eqref{eq.H2.homog.v.43} with $j = 2$. By taking derivatives to $\Lop_\eps z_\eps$ and using \eqref{eq.LePe-formula}, we have $\| \nabla \Lop_\eps z_\eps \|_2 \le C\lambda^{3/2}$. Thus the $H^3$ regularity in a bounded $C^{2,1}$ domain implies \eqref{eq.H2.homog.v.43} for $j = 3$.

\end{proof}

\begin{lemma}
    \label{l.l2.bound.z.v}
 Under the same hypotheses as in Lemma \ref{l.estimate.homog.v.43}, there are constants $C, c >0$ such that (for $v = \Lop_0^{-1}\Lop_\eps z_\eps$)
  \begin{equation}
      \label{eq.l2.bound.z.v}
      c\le (z_\eps, v) \le C 
  \end{equation}
  for all small $\eps>0$ and $z_\eps\in X_{\sqrt{\lambda},\rho}$ such that $\norm{z_\eps}_{2}=1$.
\end{lemma}

\begin{proof}
Since $z_\eps \in  X_{\sqrt{\lambda},\rho}$ is nonzero, we know that for $\rho < p_0/3$ there is $$\frac{2p_0}{3}<\sqrt{\lambda}< \sqrt{N_\ep} +\rho.$$
Now, because $z_\eps$ is a unit vector in $X_{\sqrt{\lambda},\rho}$, the proof of Lemma \ref{l.estimate.homog.v.43} already gives $\norm{\nabla z_\eps}_2 \simeq \lambda^{1/2}$. Testing the equation $\Lop_0 v = \Lop_\eps z_\eps$ by $z_\eps$, we have
\begin{equation}
    \norm{\nabla z_\eps}_2^2 \le C (\La_\eps z_\eps, z_\eps)= C(\La_0 v, z_\eps) \le C \norm{\nabla z_\eps}_2 \norm{\nabla v}_2,
\end{equation}
This gives $\norm{\nabla v}_2 \ge c\lambda^{1/2}$. Together with \eqref{eq.H2.homog.v.43}, we have
\begin{equation}
    c\lambda \le \norm{\nabla v}_2^2 \le C\lambda,
\end{equation}
for some constants $C,c>0$ independent of $\eps$ and $\lambda$. 

    Next, testing the equation $\Lop_0 v = \lambda z_\eps + r_\eps$ and using \eqref{eq.bound.reps.43} gives
    \[
    \begin{split}
          \lambda(z_\eps, v) &= (v,\Lop_0 v) - (r_\eps, v)\\
   &\in [(c-\rho) \lambda, (C +\rho)\lambda].
    \end{split}
   \]
   This proves the lemma by choosing $\rho_0 \le \min\{ c/2, p_0/3, 1/2 \}$.
\end{proof}

\subsection{Second-order correctors}

In this subsection we recall some basic properties of the first and second order correctors. We identify the 1-periodic functions as functions defined on the torus $\T^d$. Define the mean value on $\T^d$ by $\avg{f} = \fint_{\T^d} f$.
For $1\le j \le d$ the first order 1-periodic correctors $\chi_j$ satisfying
\begin{equation}
    \label{eq.corrector.1}
    -\pt_i (a_{ij} + a_{ik} \pt_k \chi_j) = 0 \ \textup{ and } \ \avg{\chi_j} = 0 ,
\end{equation}
and recall that the homogenized matrix $\hat A = (\hat a_{ij})_{d\times d}$ takes the form
\begin{equation}
    \label{eq.homogenized.tensor}
    \hat a_{ij} := \avg{a_{ij} + a_{ik} \pt_k \chi_j}.
\end{equation}

We define for $1\le i , j \le d$  
\begin{equation}
    \label{eq.define.betaij}
    \beta_{ij}: = a_{ij} + a_{ik} \pt_k \chi_j- \hat a_{ij}
\end{equation}
and a 1-periodic solution $\mu_{ij}$ that solves the following second-order corrector equation
\begin{equation}
    \label{eq.corrector.2}
    -\pt_k(a_{kl} \pt_l \mu_{ij}) = \beta_{ij} + \pt_k(a_{ki}\chi_j) \ \textup{ and }\ \avg{\mu_{ij}} = 0.
\end{equation}
The Schauder theory implies that $\chi_j$ and $\mu_{ij}$ are both $C^{1,\alpha}$ for any $\alpha \in (0,1)$.

The following lemma is classical (see \cite[Proposition 3.1.1]{shen_periodic_2018}) and we provide a proof for completeness.
\begin{lemma}
    \label{l.order.lift.tensor}
    Let $\theta_{i_1 \cdots i_k}$ be a $L^2$-integrable $k$-tensor field on $\mathbb T^d$ with $1\le i_l\le d$ for all $1\le l\le k$. Suppose 
    \begin{equation}
        \label{eq.order.lift.tensor.assume}
        \pt_{i_1} \theta_{i_1 i_2 \cdots i_k} = 0 \ \textup{ and }\ \avg{\theta_{i_1 \cdots i_k}} = 0,
    \end{equation}
    then there is a  $(k+1)$-tensor field $\omega_{ji_1 \cdots i_k } \in H^1(\mathbb T^d)$ having zero average such that 
    \begin{equation}
        \label{eq.order.lift.tensor}
        \theta_{i_1 i_2 \cdots i_k} = \pt_j \omega_{ji_1 \cdots i_k } \textup{ and }\omega_{ji_1 \cdots i_k }=-\omega_{i_1j i_2\cdots i_k }.
    \end{equation}
\end{lemma}

\begin{proof}

Since $\theta_{i_1 \cdots i_k}$ has average zero on the whole torus, we can solve for the unique solution $h_{i_1 \cdots i_k}$ (up to an additive constant) to the following equation
\begin{equation}
    \label{eq.construct.elliptic.theta}
    \Delta h_{i_1 \cdots i_k} = \theta_{i_1 \cdots i_k}.
\end{equation}
Now, we define 
\[
\omega_{j i_1 \cdots i_k} : = \partial_j h_{i_1 \cdots i_k} - \partial_{i_1} h_{j i_2\cdots i_k} \in H^1(\T^d) \textup{ for each }1\le j, i_1,\dots, i_k \le d.
\]
Note that the anti-symmetry in \eqref{eq.order.lift.tensor} follows from the definition and we also have
\[
\pt_j \omega_{ji_1 \cdots i_k } = \Delta h_{i_1 \cdots i_k} - \partial_{i_1} \partial_j h_{j i_2\cdots i_k} = \theta_{i_1 \cdots i_k},
\]
where the term $\partial_{i_1} \partial_j h_{j i_2\cdots i_k}$ vanishes because by \eqref{eq.construct.elliptic.theta} 
\[
 \Delta \partial_jh_{j i_2 \cdots i_k} =\partial_j \theta_{j i_2 \cdots i_k} =0,
\]
which shows that $\partial_jh_{j i_2 \cdots i_k}$ must be a constant.
    
\end{proof}

By Lemma~\ref{l.order.lift.tensor}, there is an $H^1$ 3-tensor field $\gamma_{kij}$ for $1\le k,i,j\le d$ such that 
\begin{equation}
    \label{eq.define.gamma.kij}
   \beta_{ij} = \pt_k \gamma_{kij} \textup{ and }\gamma_{kij} = -\gamma_{ikj}.
\end{equation}
The $C^\alpha$ regularity of $\beta_{ij}$ implies that $\gamma_{kij}$ is $C^{1,\alpha}$, due to the proof of Lemma \ref{l.order.lift.tensor}.
Using $\gamma_{kij}$, we can rewrite the second order corrector equation \eqref{eq.corrector.2} as follows
\begin{equation}
    \label{eq.corrector.2.gamma}
    \partial_k ~\underbrace{\lb a_{kl}\pt_l \mu_{ij} + \gamma_{kij} + a_{ki}\chi_j  \rb}_{\xi_{kij}^*} = 0.
\end{equation}

We denote 
\begin{equation}\label{eq.define.xi.kij}
    \xi_{kij} : = \xi_{kij}^* - \avg{\xi_{kij}^*}
\end{equation}
and by Lemma~\ref{l.order.lift.tensor} again, we obtain a 4-tensor $\zeta_{lkij}$ on the torus such that 
\begin{equation}
    \label{eq.4tensor.zeta}
    \xi_{kij} = \partial_l \zeta_{lkij} \textup{ and } \zeta_{lkij} = - \zeta_{klij}.
\end{equation}

\begin{lemma}
    \label{l.cancellation.3.tensor}
    The 3-tensor $\avg{\xi_{kij}^*}$ satisfies
    \begin{equation}
        \label{eq.symm.zero.3}
       \avg{\xi_{kij}^*} = - \avg{\xi_{jik}^*} 
    \end{equation}
\end{lemma}

\begin{proof}
This property is known (see e.g. \cite[p.~1271]{MV97}) and we provide a proof for the reader's convenience.
Testing the first-order corrector equation  \eqref{eq.corrector.1} against the second-order correctors $\mu_{ij}$, we obtain
\begin{equation}
    \label{eq.cancellation.3.1}
    \avg{a_{kl}\pt_l \mu_{ij}} = -\avg{ \mu_{ij} \La_1 \chi_k} = -\avg{\chi_k \La_1 \mu_{ij}}.
\end{equation}
By the second corrector equation \eqref{eq.corrector.2}, we obtain
\begin{equation}
    \label{eq.cancellation.3.2}
    \begin{split}
       \avg{a_{kl}\pt_l \mu_{ij}} &= -\avg{\chi_k \lb \beta_{ij} + \pt_l(a_{li}\chi_j)\rb}\\
  &=-\avg{\chi_k\lb a_{ij} + a_{il} \pt_l \chi_j- \hat a_{ij} \rb} + \avg{\pt_l \chi_k a_{li} \chi_j}\\
  &=-\avg{\chi_k a_{ij} + \chi_k a_{il} \partial_l \chi_j - \partial_l \chi_k a_{li} \chi_j },
    \end{split}
   \end{equation} 
where we have used the fact that $\avg{\chi_k} = 0$. Since $\avg{\gamma_{kij}} = 0$, we obtain
\[
\avg{\xi_{kij}^*} = \avg{a_{ki}\chi_j - \chi_k a_{ij}} + \avg{ \partial_l \chi_k a_{li} \chi_j - \chi_k a_{il} \partial_l \chi_j},
\] 
which proves the lemma.
\end{proof}

\subsection{Second-order approximation scheme}

For any 1-periodic function $h$ (or equivalently any function on the torus), we denote 
\begin{equation}
    \label{eq.def.eps.scale}
    h^\eps(x) : = h\lb\frac{x}{\eps}\rb.
\end{equation}

Suppose $v\in H^3\cap H_0^1(\Om)$, we denote $z_\eps = \La_\eps^{-1}\La_0 v $ (solved with zero boundary condition), and start to work with the following approximation ansatz
\begin{equation}
    \label{eq.second.order.ansatz}
z_\eps = v + \eps \chi_j^\eps \pt_j v + \eps B_{1,\eps} + \eps^2 \mu_{ij}^\eps\pt_{ij}v + \eps^2 B_{2,\eps} + R_{2,\eps},
\end{equation}
where the $B_{i,\eps}$'s are boundary layer terms defined as follows respectively
\begin{equation}
    \label{eq.boundary.layers.def}
    B_{1,\eps} : = \cH_\eps \lmb - \chi_j^\eps \restr{\pt_jv}{\pt\Om}\rmb \textup{ and }B_{2,\eps} : = \cH_\eps\lmb-\mu_{ij}^\eps \restr{\pt_{ij}v}{\pt\Om}\rmb.
\end{equation}
Here $\cH_\eps$ denotes the $\Lop_\eps$-harmonic extension of the input boundary data in the whole domain $\Om$; see \eqref{eq:Poisson-integral} for a precise formula in terms of Poisson integral.

\begin{lemma}
    \label{l.second.error}
    Let $v\in H^3\cap H_0^1(\Om)$ and $z_\eps = \La_\eps^{-1}\La_0 v $ with zero boundary condition. Then the remainder $R_{2,\eps} \in H_0^1(\Omega)$ in \eqref{eq.second.order.ansatz} has the following estimate
    \begin{equation}
        \label{eq.second.error}
        \norm{R_{2,\eps}}_{H_0^1(\Om)} \le C \eps^2 \norm{v}_{H^3(\Om)},
    \end{equation}
    where the constant $C>0$ is independent of $\eps$.
\end{lemma}

\begin{proof}
    This follows from Lemma \ref{l.cancellation.3.tensor} and a similar calculation as \cite[Theorem 9.1]{SZ18} or \cite[Theorem 5.1]{Z20}.
    The layers $B_{i,\eps}$ as defined in \eqref{eq.boundary.layers.def} cancel the corrections on the boundary, and then $R_{2,\eps}$ belongs to $H_0^1(\Omega)$. To prove the estimate \eqref{eq.second.error}, we observe that 
\begin{equation}
    \label{eq.equation.R.2.eps}
    \begin{split}
        \La_\eps R_{2,\eps} & = \La_\eps \lb z_\eps - v -  \eps \chi_j^\eps \pt_j v - \eps^2 \mu_{ij}^\eps\pt_{ij}v \rb\\
        & = (\La_0-\La_\eps) v + \pt_i(a_{ik}^\eps(\pt_k\chi_j)^\eps \pt_j v) +\eps\pt_i(a_{ik}^\eps \chi_j^\eps \pt_{jk}v) +\eps^2 \La_\eps(\mu_{ij}^\eps\pt_{ij}v)\\
        &=\pt_i(\beta_{ij}^\eps \pt_j v)+\eps\pt_i(a_{ik}^\eps \chi_j^\eps \pt_{jk}v) +\eps^2 \La_\eps(\mu_{ij}^\eps\pt_{ij}v).
    \end{split}
\end{equation}
Note that by the anti-symmetry of $\gamma$ in \eqref{eq.order.lift.tensor}
\[
   \pt_i(\beta_{ij}^\eps \pt_j v) = \beta_{ij}^\eps \pt_{ij}v=\eps \pt_k \gamma_{kij}^\eps \pt_{ij} v=\eps \pt_k(\gamma_{kij}^\eps \pt_{ij} v).
\]
This, combining \eqref{eq.equation.R.2.eps}, the anti-symmetry of $\avg{\xi_{kij}^*}$ in Lemma~\ref{l.cancellation.3.tensor}, and the same property of $\zeta_{lkij}$ as defined in \eqref{eq.4tensor.zeta}, shows that
\begin{equation}
    \label{eq.equation.R.2.eps.ctd}
    \begin{split}
        \La_\eps R_{2,\eps} & = -\eps \pt_k\lb\lb\xi_{kij}^*\rb \pt_{ij} v\rb - \eps^2 \pt_l(a_{lk}^\eps \mu_{ij}^\eps \pt_{kij} v)\\
        &= -\eps \xi_{kij}^* \pt_{kij} v - \eps^2 \pt_l(a_{lk}^\eps \mu_{ij}^\eps \pt_{kij} v)\\
        &= -\eps \xi_{kij} \pt_{kij} v- \eps^2 \pt_l(a_{lk}^\eps \mu_{ij}^\eps \pt_{kij} v)\\
        &= -\eps^2 \pt_l \lb\lb\zeta_{lkij}^\eps + a_{lk}^\eps \mu_{ij}^\eps\rb \pt_{kij} v\rb.
    \end{split}
\end{equation}
This proves the lemma by standard $H^1$ estimates and the $L^\infty$ boundedness of $\zeta$ and $\mu$.
\end{proof}

\begin{remark}
    \label{r.barrier.of.higher.order.expansion}
    We remark here a major barrier of utilizing higher order homogenization expansions. A heuristic calculation shows that if one use a $k$-th order analogous expansion to \eqref{eq.second.order.ansatz} (with smooth coefficients and domain), the formal $k$-th error $R_{k,\eps}$ should satisfy
    \[
      \norm{R_{k,\eps}}_{H_0^1(\Om)} \le C \eps^k \norm{v}_{H^{k+1}(\Om)}.
    \]
    Unfortunately, the function $v$, solving $\La_0 v = \La_\ep z_\eps \in X_{\sqrt{\lambda},\rho}$, only has the following bound (in the full subcritical regime $\lambda \le C \eps^{-2}$)
    \[
    \norm{v}_{H^{k+1}(\Om)} \le C_k \lambda \norm{z_\eps}_{H^{k-1}(\Om)} \le C_k \ep^{2-k} \lambda^{3/2}.
    \]
    This only gives for $k\ge 2$, $\norm{R_{k,\eps}}_{H_0^1(\Om)} \le C_k \eps^2\lambda^{3/2} $, which is not improving the second order regime.
\end{remark}

\begin{lemma}
    \label{l.boundarylayer.in.v}
    The boundary layer terms $B_{i,\eps}$ in \eqref{eq.second.order.ansatz} satisfy the following estimates
    \begin{equation}
        \label{eq.boundarylayer.in.v}
        \norm{B_{i,\eps}}_{H^1(\Om)} \le C \lb \eps^{-1/2} \lb\norm{v}_{H^{i}(\Om)} \norm{v}_{H^{i+1}(\Om)}\rb^{1/2} + \norm{v}_{H^{i+1}(\Om)} \rb,
    \end{equation}
    for $i = 1, 2$.     For the $L^2$ norm, the estimate can be improved to
     \begin{equation}
        \label{eq.boundarylayer.in.v.L2}
        \norm{B_{1,\eps}}_{L^2(\Om)} \le C \lb\norm{v}_{H^1(\Om)}\norm{v}_{H^{2}(\Om)}\rb^{1/2} \ \textup{ and }\ \norm{B_{2,\eps}}_{L^2(\Om)} \le C\norm{v}_{H^3(\Om)}.
    \end{equation}
    Here the constants $C>0$ are independent of $\eps$.
\end{lemma}

\begin{proof}

The $H^1$ estimate \eqref{eq.boundarylayer.in.v} is obtained by applying Lemma~\ref{l.multiplicative.trace.theorem.apx} and Lemma~\ref{l.H1half.boundary.trace.apx}. The proof of \eqref{eq.boundarylayer.in.v.L2} follows from \cite[Theorem 3]{AvellanedaLin} and Lemma~\ref{l.multiplicative.trace.theorem.apx} and the standard trace theorem.

\end{proof}

By combining Lemma~\ref{l.second.error} and Lemma~\ref{l.boundarylayer.in.v}, we obtain the following $L^2$ estimate.

\begin{corollary}
    \label{cor.L2.difference.z.v}
    The difference $z_\ep-v$ satisfies the following estimate
    \begin{equation}
        \label{eq.L2.difference.z.v}
        \norm{z_\ep-v}_{L^2(\Om)} \le C \lb \eps \norm{v}_{H^1} + \eps \lb\norm{v}_{H^1}\norm{v}_{H^2}\rb^{1/2} + \eps^2 \norm{v}_{H^3}\rb.
    \end{equation}
    Here the constant $C$ is independent of $\eps$, $z$ and $v$.
\end{corollary}

Here and after, for $t>0$, write
\[
\Om_t : = \{x\in \Om \;:\; \dist(x,\pt \Om) < t\}.
\]

\begin{lemma}
    \label{l.grad.error.z.v}
Assume the same hypotheses as Lemma \ref{l.second.error}. 
There is a sequence of matrix fields $J_\eps$ such that for any $b > 0$,
\begin{equation}
    \label{eq.grad.error.z.v}
    \norm{\gd z_\eps  - J_\eps \gd v}_{L^2(\Om_{b\eps})} \le C_b \lb \eps^{3/2}\lb\norm{v}_{H^2(\Om)}\norm{v}_{H^3(\Om)}\rb^{1/2} + \ep^2 \norm{v}_{H^3(\Om)} \rb,
\end{equation}
Moreover, there exists $b_0 > 0$ such that the matrix field $J_\eps$ satisfies $ \norm{J_\eps}_{L^\infty(\Om)}\le C$ and for some constant $c>0$
\begin{equation}
    \label{eq.J.eps}
  |J_{\ep}(x)\xi| \ge c |\xi| \quad \textup{ for all } x\in \Om_{b_0 \ep}, \ 0<\ep<\ep_0 \textup{ and }\xi\in\R^d.
\end{equation}
The constants $C, c, b_0, \eps_0$ are all independent of $\eps$ and $z_\eps$.

\end{lemma}

\begin{proof}
Consider the following equation for the Dirichlet correctors $\Phi_\eps: = (\Phi_{j,\ep})_{j=1}^d$
\begin{equation}
    \label{eq.Dirichlet.correct}
    \La_\ep \Phi_{j,\eps} = 0 \quad \textup{in }\Om \quad \textup{and} \quad \Phi_{j,\eps} = x_j \quad \textup{on }\pt\Om.
\end{equation}
We define 
\begin{equation}
    \label{eq.define.Jeps}
    J_{\eps}  : = D \Phi_\eps = (\partial_i \Phi_{j,\eps})_{d\times d}.
\end{equation}
The $L^\infty(\Om)$-boundedness of $J_\ep$ and the coercivity property in \eqref{eq.J.eps} follow from \cite[Lemma 5.1 and (5.5)]{kenig_estimates_2013}.

To prove \eqref{eq.grad.error.z.v}, we write $\Theta_{j,\eps}$ as the $\La_\ep$-harmonic extension of $-\restr{\chi_j^\eps}{\pt\Om}$, i.e., $\Theta_{j, \eps} = - \cH_\eps\lb\restr{\chi_j^\eps}{\pt\Om}\rb$. Observe that 
\[
\Phi_{j,\eps} = x_j + \ep \chi_j^\ep + \ep \Theta_{j,\eps}.
\]
Writing $W_\ep = B_{1,\eps} - \Theta_{j,\eps} \pt_j v$, we obtain
\[
v + \ep \chi_j^\ep \pt_j v + \ep B_{1,\ep} = v + (\Phi_{j,\ep} - x_j) \pt_j v + \ep W_\eps.
\]
   Observe that
   \begin{equation}
       W_\eps = \cH_\eps (\restr{-\chi_j^\eps \partial_j v}{\partial \Om}) + \partial_j v \cH_\eps (\restr{\chi_j^\eps}{\partial \Om}) = - [\cH_\eps, \partial_j v] \chi_j^\eps.
   \end{equation}

Note that the gradient of the above display takes the following form
\begin{equation}
    \label{eq.grad.first.scheme}
    \gd \lb v + \ep \chi_j^\ep \pt_j v + \ep B_{1,\ep}\rb = J_\eps \gd v + (\Phi_{j,\ep} - x_j) \pt_j \gd v + \ep \gd W_\eps.
\end{equation}
We now apply the second order expansion \eqref{eq.second.order.ansatz} to $z_\eps$ and obtain
\begin{equation}
    \label{eq.gd.z.v.boundary.collar.estimate}
    \begin{split}
        \gd z_\ep - J_\ep \gd v &=\underbrace{(\Phi_{j,\ep} - x_j) \pt_j \gd v}_{H_1} ~+~ \underbrace{\ep \gd W_\eps}_{H_2}~+~ \underbrace{\ep^2\gd (\mu_{ij}^\ep \pt_{ij}v)}_{H_3}\\[10pt]
        & + \underbrace{\ep^2 \gd B_{2,\eps}}_{H_4} ~+~ \underbrace{\gd R_{2,\eps} }_{H_5}.
    \end{split}
   \end{equation}

By combining the commutator estimate in Theorem \ref{thm:commutator}, Lemma \ref{l.multiplicative.trace.theorem.apx} and Lemma \ref{l.boundary.collar.gradient.estimate.apx}, there is 
\begin{equation}
    \label{eq.bound.H2}
    \begin{aligned}
     \int_{\Om_{b\ep}} |H_2|^2 \dd x & \le C \eps^3 \norm{\nabla v}_{L^2(\partial \Omega)}^2 + C\e^2 \norm{\nabla^2 v}_{L^2(\Omega_{2b\ep})}^2 \\
     & \le C\eps^3 \norm{v}_{H^2(\Omega)} \norm{v}_{H^3(\Omega)}.
     \end{aligned}
\end{equation}

For $H_1$ and $H_3$ we observe that by Lemma \ref{l.boundary.collar.gradient.estimate.apx} 
\begin{equation}
    \label{eq.boundary.collar.second.grad}
    \int_{\Om_{b\ep}} |\nabla^2 v|^2 \dd x \le C \ep \norm{v}_{H^2(\Om)}\norm{v}_{H^3(\Om)},
\end{equation}
which shows that (note that $\sup_j\norm{\Phi_{j,\ep}-x_j}_{L^\infty(\Omega)}\le C \eps$)
\begin{equation}
    \label{eq.bound.H1.H3}
    \begin{split}
         \int_{\Om_{b\ep}} |H_1|^2+|H_3|^2 \dd x &\le C\ep^2 \lb \int_{\Om_{b\ep}} |\nabla^2 v|^2 \dd x + \ep^2\norm{v}_{H^3(\Om)}^2 \rb\\
&\le C  \ep^3 \norm{v}_{H^2(\Om)}\norm{v}_{H^3(\Om)} + C\ep^4\norm{v}_{H^3(\Om)}^2.
    \end{split}
    \end{equation}

Lastly, by Lemma~\ref{l.second.error} and Lemma~\ref{l.boundarylayer.in.v} 
\begin{equation}
    \label{eq.bound.H4.H5}
    \int_{\Om_{b\ep}} |H_4|^2+|H_5|^2 \dd x \le C  \ep^3 \norm{v}_{H^2(\Om)}\norm{v}_{H^3(\Om)} + C\ep^4\norm{v}_{H^3(\Om)}^2.
\end{equation}
This finishes the proof of the lemma.
\end{proof}

\subsection{Boundary estimates of \texorpdfstring{$v$}{v}}
\label{sec.Dv-est}
In this subsection, we prove the upper and lower boundary gradient estimates of $v = \Lop_0^{-1} \Lop_\eps z_\eps$. The following lemma is the main result of this subsection.

\begin{lemma}
    \label{l.estimate.homog.v.boundary.43}
 Suppose $\Omega$ is $C^{2,1}$, $z_\eps\in X_{\sqrt{\lambda},\rho}$ is a unit vector and $\rho < \rho_0$ for some small $\rho_0>0$. Let $v = \Lop^{-1}_0 \Lop_\eps z_\eps$. For $N_\eps \le \delta\eps^{-4/3}$ with sufficiently small $\delta>0$, we can find constants $c, C>0$ such that 
 \begin{equation}
     \label{eq.boundaryestimate.homog.v}
      c \sqrt{\lambda} \le \norm{B v}_{Y} = \norm{\gd v}_{L^2(\pt \Omega;\mathbb C)} \le  C \sqrt{\lambda}. 
 \end{equation}

\end{lemma}

The next corollary converts the above boundary estimate to an estimate in the thin layer. This allows us to use the comparison in Lemma \ref{l.grad.error.z.v}.
\begin{corollary}
    \label{cor.collar.estimate.v}
    Assume the same hypotheses as in Lemma~\ref{l.estimate.homog.v.boundary.43}. Then for all $a>0$ the function $v$ satisfies 
    \begin{equation}
        \label{eq.collar.estimate.v}
        c {\lambda} \le \frac{1}{ \eps}\int_{\Om_{a\eps}} |\gd v|^2 \dd x \le C {\lambda},
    \end{equation}
    for some constants $C,c >0$ independent of $\eps$ (depending on $a$).
\end{corollary}

\begin{proof}
  
    By Lemma \ref{lem.layer2bdry}, we have
    \begin{equation}
        \label{eq.modified.trace1}
         \frac{1}{a\eps} \int_{\Om_{a\eps}} |\gd v|^2 \dd x \le C \int_{\pt\Om} |\gd v|^2 \dd \sigma + C a\eps \norm{\nabla^2v}_{L^2(\Om_{a\eps})}^2.
    \end{equation}
    \begin{equation}
        \label{eq.modified.trace2}
        \frac{1}{a\eps} \int_{\Om_{a\eps}} |\gd v|^2 \dd x \ge c \int_{\pt\Om} |\gd v|^2 \dd \sigma - C a\eps \norm{\nabla^2v}_{L^2(\Om_{a\eps})}^2.
    \end{equation}
    By Lemma~\ref{l.boundary.collar.gradient.estimate.apx} and Lemma~\ref{l.estimate.homog.v.43}
\[
\norm{\nabla^2v}_{L^2(\Om_{a\eps})}^2 \le C a\eps\norm{v}_{H^2(\Om)}\norm{v}_{H^3(\Om)} \le C a\eps\lambda^{1+3/2}.
\]
This proves the lemma by combining \eqref{eq.modified.trace1}, \eqref{eq.modified.trace2} and Lemma~\ref{l.estimate.homog.v.boundary.43}. 
\end{proof}

Let us begin with the following auxiliary lemma.

\begin{lemma}
    \label{l.estimate.L2.size.v.43}
    Assume the same hypotheses as in Lemma~\ref{l.estimate.homog.v.boundary.43}. The function $v$ satisfies
    \begin{equation}
        \label{eq.estimate.L2.size.v.43}
        \frac12 \le \norm{v}_{L^2(\Omega)} \le 2.
    \end{equation}
\end{lemma}

\begin{proof}
    By Lemma~\ref{cor.L2.difference.z.v}, we know that 
    \[
    \norm{z_\ep-v}_2 \le C \lb \ep N_\ep^{1/2} + \ep N_\ep^{3/4} + \ep^2 N_{\ep}^{3/2} \rb \le C \delta^{3/4}
    \]
    when $\delta>0$ is chosen sufficiently small. Since $z_\eps$ is a unit vector, the lemma is thus proved.
\end{proof}

\begin{proof}[Proof of Lemma~\ref{l.estimate.homog.v.boundary.43}]

We recall the notation in \eqref{eq.homog.v.43} and write
\[
\La_0 v = \La_\ep z_\eps = \lambda z_\eps + r_\eps.
\]
By the Rellich identity, we obtain for any smooth vector field $h$ there is 
\begin{equation}
    \label{eq.rellich.general.h}
    \begin{split}
        \int_{\pt \Omega}
(h\cdot \nu)\,
\gd v\cdot \widehat{A}\gd v
\dd \sigma
&=
2\left(\widehat{A}\gd v,\, \gd h\,\gd v\right)
-\left(\operatorname{div}h\,\widehat{A}\gd v,\, \gd v\right)
\\
&\quad-2\lambda\left(h\cdot \gd v,\, z_\eps\right)
-2\left(h\cdot \gd v,\, r_\eps\right).
    \end{split}
\end{equation}
To prove the upper bound in \eqref{eq.boundaryestimate.homog.v}, we observe that since $\Omega$ is $C^{2,1}$, there is a $C^{1,1}$ vector field $h$ satisfying $(h\cdot \nu) \equiv 1 >0$ uniformly on the boundary of the domain. Therefore 
\begin{equation}
    \label{eq.rellich.upper.bound}
    \begin{split}
         \norm{\gd  v}_{L^2(\pt\Om)}^2 &\le C \int_{\pt \Omega}
(h\cdot \nu)\,
\gd v\cdot \widehat{A}\gd v
\dd \sigma \\
& \le C \lambda + 2\lambda ~ \underbrace{|\left(h\cdot \gd v,\, z_\eps\right)|}_{K_1}.
    \end{split}
   \end{equation}
Here in the last line we have used \eqref{eq.bound.reps.43} to bound the terms involving $r_\ep$. Now we observe that by Lemma \ref{l.estimate.L2.size.v.43}
\begin{equation}
    \label{eq.bound.K1.upperbound}
    \begin{split}
          |K_1| & \le |(h\cdot \gd v, v)| + |(h\cdot \gd v, z_\ep - v)|\\
&\le C \norm{v}_{2}^2 + |(h\cdot \gd v, z_\ep - v)|\\
&\le C + |(h\cdot \gd v, z_\ep - v)|.
    \end{split}
  \end{equation}
By using the expansion \eqref{eq.second.order.ansatz}
\begin{equation}
    \label{eq.bound.J3.decomp}
        \begin{split}
|(h\cdot \gd v, z_\ep - v)| &\le  \ep |(h\cdot \gd v, \chi_j^\ep \pt_j v)| + \ep|(h\cdot \gd v, B_{1,\ep}) |\\
& \quad + \ep^2|(h\cdot \gd v, \mu_{ij}^\ep \pt_{ij}v)| + \ep^2|(h\cdot \gd v, B_{2,\ep})| \\
&\quad + |(h\cdot \gd v, R_{2,\eps}) |\\
& =: K_{11}+K_{12}+K_{13}+K_{14} +K_{15} .
    \end{split}
\end{equation}
For $K_{11}$, we know that $\chi_j$ have zero average, and then there are $\ep$-periodic fields $\Psi_{ij}^\ep$ such that $\chi_j^\ep = \ep \pt_i \Psi_{ij}^\ep$. This shows that by Lemma~\ref{l.estimate.homog.v.43}
\[
\begin{split}
    K_{11}& \le \ep^2|(\pt_i(h\cdot \gd v\pt_j v), \Psi_{ij}^\ep)| + \ep^2 \lw\int_{\pt\Om}\nu_i\Psi_{ij}^\ep \pt_j v h\cdot \gd v \dd \sigma\rw\\
    &\le C \ep^2(\lambda^{3/2}+ \lambda) + C\ep^2  \norm{\gd  v}_{L^2(\pt\Om)}^2 .
\end{split}
\]
As $\lambda \le CN_\ep < C\delta \eps^{-4/3}$, we have
\begin{equation}
    \label{eq.bound.K11}
    K_{11} \le C \delta^{3/2} + C \eps^2 \norm{\gd  v}_{L^2(\pt\Om)}^2 .
\end{equation}
For $K_{12}$, we have by Lemma~\ref{l.estimate.L2.size.v.43} and Lemma~\ref{l.fine.estimate.rellich.integral.apx}
\[
K_{12} \le C\ep \norm{\gd v}_{L^2(\pt\Om)} \lb\norm{v}_{H^1(\Om)}^{1/2} + C\rb \le C \ep \lambda^{1/4} \norm{\gd v}_{L^2(\pt\Om)}.
\]
For $K_{13}, K_{14}$ and $K_{15}$, we have by Lemma~\ref{l.second.error}, Lemma~\ref{l.fine.estimate.rellich.integral.apx} and $\lambda < C\delta \ep^{-4/3}$
\[
\begin{split}
    K_{13}+K_{14}+K_{15}&\le \eps^2 \lb\norm{v}_{H^1}\norm{v}_{H^2} + \norm{B_{2,\ep}}_{L^2(\pt\Om)}\lb\norm{v}_{H^1(\Om)}^{1/2} +C\rb\rb\\
    &\quad + \norm{R_{2,\ep}}_{H^1}\\
&\le C \eps^2 (\lambda^{3/2} + \lambda^{5/4})\\
&\le C \delta^{3/2}.
\end{split}
\]
Combining \eqref{eq.rellich.upper.bound}, \eqref{eq.bound.K1.upperbound}, \eqref{eq.bound.J3.decomp} and the estimates of $K_{1j}$, we obtain if we denote $M : = \norm{\gd v}_{L^2(\pt\Om)}$,
\[
\begin{split}
    M^2 &\le C \lambda  + C \lambda \lb \ep^2 M^2 + \ep \lambda^{1/4} M + \delta^{3/2}\rb\\
    &\le C \lambda + C \lambda( (\ep^{4/3}+\ep^2)M^2 +C\delta^{1/2}),
\end{split}
\]
which proves that $M \le C \sqrt{\lambda}$ as $\lambda(\ep^{4/3}+\ep^2) \le C\delta \ll 1$.

To prove the lower bound in \eqref{eq.boundaryestimate.homog.v}, we choose $h (x) = x$ in \eqref{eq.rellich.general.h} and obtain the following lower bound for $M$
        \begin{equation}
        \label{eq.Rellich.identity.alt}
        \begin{split}
            \int_{\pt \Omega} (x\cdot \nu)~ \gd v \cdot \widehat{A} \gd v \dd \sigma & = (2-d)(\gd v , \widehat{A }\gd v) -2 \lambda (x\cdot \gd v , z_\eps) - 2(x\cdot \gd v, r_\eps)\\
            &=(2-d)(v,\Lop_0 v)-2 \lambda (x\cdot \gd v , z_\eps)- 2(x\cdot \gd v, r_\eps)\\[5pt]
            &= 2 \lambda(v,z_\eps)  -d \lambda (v,z_\eps - v) + (2-d) (v,r_\eps) \\[5pt]
            &\quad- 2\lambda (x\cdot \gd v, z_\eps -v)- 2(x\cdot \gd v, r_\eps)\\[5pt]
            &=:2 \lambda(v,z_\eps) + J_1 + J_2 + J_3 +J_4.
        \end{split}
    \end{equation}

It is not difficult to derive by \eqref{eq.bound.reps.43} and \eqref{eq.H2.homog.v.43} 
\begin{equation}
    \label{eq.bound.easy.J}
    |J_2|+|J_4|\le C_d \norm{\gd v}_{2} \norm{r_\eps}_{2} \le C \rho \lambda.
\end{equation}
By Lemma~\ref{l.l2.bound.z.v}, both $J_2$ and $J_4$ can be absorbed into the leading term
\[
2 \lambda(v,z_\eps) \simeq \lambda
\]
for sufficiently small bandwidth $\rho$.

It then suffices to bound $J_1$ and $J_3$. For $J_1$ we have by combining Lemma~\ref{l.estimate.homog.v.43}, Corollary~\ref{cor.L2.difference.z.v} and Lemma~\ref{l.estimate.L2.size.v.43}
\begin{equation}
    \label{eq.bound.J1}
          |J_1|  \le d~ \lambda  \norm{v}_2 \norm{ z_\ep - v}_2   \le C \delta^{3/4} \lambda .
  \end{equation}
On the other hand, the bound $|J_3| \le C \delta^{3/4}\lambda$ is obtained similarly as \eqref{eq.bound.J3.decomp} with the additional upper bound $\norm{\gd  v}_{L^2(\pt \Om)} \le C \sqrt{\lambda}$ that we have just obtained. This finishes the proof of the lower bound once one chooses small $\delta>0$.

\end{proof}

\subsection{Boundary gradient estimates of \texorpdfstring{$z_\ep$}{z} and proofs of main theorems}

In this subsection, we prove the upper and lower bounds of the boundary gradient $\int_{\pt \Om} |\gd z_\eps|^2 \dd \sigma$, followed by the proof of Theorem \ref{intro:cluster} (and therefore Theorem \ref{intro:main}).

We begin with a two-sided boundary layer estimate for $\nabla z_\eps$.
\begin{lemma}
    \label{l.boundary.estimate.Z.ep}
    Assume the same hypotheses as in Theorem \ref{intro:cluster} and $\| z_\eps \|_2 = 1$. Then for all $a>0$ there are positive constants $C, c >0$ (depending on $a$) such that for $0<\ep<\ep_0$
\begin{equation}
    \label{eq.boundary.estimate.Z.ep.1}
    c{\lambda}\le \frac{1}{\eps}\int_{\Om_{a\eps}} |\gd z_\eps|^2 \dd x \le C {\lambda}.
\end{equation}
Here $\Om_t:=\{x\in\Om\;:\;\dist(x,\pt\Om)<t\}$. The constants $C, c$ are independent of $\eps$, $\lambda$ and $z_\eps\in X_{\sqrt{\lambda},\rho}$.

\end{lemma}

\begin{proof}
    The proof is finished by combining Lemma~\ref{l.estimate.homog.v.43}, Lemma~\ref{l.grad.error.z.v} and Corollary~\ref{cor.collar.estimate.v}.
\end{proof}

\begin{lemma}
    \label{l.bl.blow.up}
Assume the same hypotheses as in Lemma \ref{l.boundary.estimate.Z.ep}. Then there is a constant $C>0$ such that
\begin{equation}
    \label{eq.boundary.control.collar}
    \int_{\pt\Om}|\gd z_\ep|^2 \dd \sigma \le C\lb\frac{1}{\ep}\int_{\Om_{2\ep}}|\gd z_\eps|^2 \dd x + \ep \int_{\Om_{2\eps}}|\La_\eps z_\eps|^2 \dd x\rb.
\end{equation}

\end{lemma}

\begin{proof}

Let $\eta_\eps$ be a smooth cutoff function such that $\eta_\eps = 1$ in $\Omega_{\eps/2}, \eta_\eps = 0$ in $\Omega \setminus \Omega_{\eps}$, and $|\nabla \eta_\eps| \le C\eps^{-1}$. Let $h$ be the same vector field as in the proof of Lemma \ref{l.estimate.homog.v.boundary.43} such that $h\cdot \nu \equiv 1$ on $\partial \Omega$. Then by the divergence theorem
\begin{equation}
\begin{aligned}\label{eq.ze-Rellich-1}
    \int_{\partial \Omega} |\nabla z_\eps|^2 \dd \sigma & = \int_{\partial \Omega} \eta_\eps h\cdot v |\nabla z_\eps|^2 \dd \sigma = \int_{\Omega} \diver (\eta_\eps h |\nabla z_\eps|^2) \\
    & \le \frac{C}{\eps} \int_{\Omega_\eps} |\nabla z_\eps|^2 + C\eps \int_{\Omega_\eps} |\nabla^2 z_\eps|^2.
\end{aligned}
\end{equation}
Now we claim that
\begin{equation}\label{eq.D2z.layer}
    \int_{\Omega_\eps} |\nabla^2 z_\eps|^2 \le C\int_{\Omega_{2\eps}} |\Lop_\eps z_\eps|^2 + \frac{C}{\eps^2} \int_{\Omega_{2\eps}} |\nabla z_\eps|^2,
\end{equation}
which together with \eqref{eq.ze-Rellich-1} gives \eqref{eq.boundary.control.collar}.

To show the claim, we first establish estimates in boundary patches.
For each $p \in \pt \Om$, let $D_{\theta}(p): = B_{\theta\eps}(p)$ for any $\theta>0$. When $\eps>0$ is sufficiently small, the domain $D_\theta(p) \cap \Om$ is Lipschitz and piecewise $C^{2,1}$ as $\pt \Om$ is assumed to be $C^{2,1}$. Without loss of generality, we assume $\theta =2$. By rescaling the boundary patch $D_2(p)$ and using the $H^2$ estimate for the operator $\La_1$ with $C^{1,1}$ boundary, 
\begin{equation}\label{eq.D2z-patch}
    \int_{D_{3/2}(p) \cap \Om}|\nabla^2 z_\eps|^2  \dd x \le C  \int_{D_{ 2}(p) \cap \Om}|\La_\eps z_\eps|^2  \dd x + \frac{C}{\eps^2}\int_{D_{ 2}(p) \cap \Om}|\gd z_\eps|^2  \dd x .
\end{equation}

Now we choose a maximally $\frac12\eps$-separated point set $S_\eps$ of $\pt \Om$ (i.e. all two distinct points in $S_\eps$ have distance at least $\frac12 \eps$, and adding any additional point from $\pt\Om$ breaks the property), and for each $p\in S_\eps$, we associate a ball $D_2(p)$ of radius $2\eps$. Note that $\{D_{1/2}(p) \}_{p\in S_\eps}$ covers $\pt \Om$ and $\{D_2(p) \cap \overline{\Om}\}_{p\in S_\eps}$ are contained in $\Om_{2\eps}$ for  all small $\eps$, and the multiplicity of the cover is bounded by a fixed number independent of $\eps$ and $S_\eps$ (due to the $C^{2,1}$-smoothness of the boundary). Moreover,
\begin{equation}
    \Omega_{\eps} \subset \bigcup_{p\in S_\eps} D_{3/2}(p) \cap \Omega.
\end{equation}
Consequently, summing \eqref{eq.D2z-patch} over all $p\in S_\eps$ proves \eqref{eq.D2z.layer}.
\end{proof}

One might wonder if there is a lower bound of the form
\[
 \int_{\pt\Om}|\gd z_\ep|^2 \dd \sigma \ge   \frac{1}{C\ep}\int_{\Om_{b\ep}}|\gd z_\eps|^2 \dd x - C\ep \int_{\Om_{a\eps}}|\La_\eps z_\eps|^2 \dd x,
\]
for some $a>b>0$. This does not immediately follow the same argument as in the proof in Lemma \ref{l.bl.blow.up}. Specifically, one can never improve \eqref{eq.D2z-patch} to
\[
 \int_{D_{3/2}(p) \cap \Om}|\nabla^2 z_\eps|^2  \dd x \le C_\delta  \int_{D_{ 2}(p) \cap \Om}|\La_\eps z_\eps|^2  \dd x - \frac{\delta}{\eps^2}\int_{D_{ 2}(p) \cap \Om}|\gd z_\eps|^2  \dd x 
\]
for arbitrary small $\delta>0$, which is obviously false if $z_\eps$ is a nontrivial $\Lop_\eps$-harmonic function.

To compensate this, one needs to consider a local boundary energy coercivity principle.

Let $Z_s = B_s(p) \cap \Omega$ for some $p \in \partial \Omega$, and $I_s = B_s(p) \cap \partial \Omega$.
\begin{lemma}
    \label{l.local.coercivity}
    Fix $\theta>0$. Suppose \begin{equation}\label{eq.local.coercivity.assume}
        -\textup{div}(A \gd u) = F \textup{ in }Z_2, \quad \  u=0 \textup{ on }I_2, \quad \ F \in L^2(Z_2).
    \end{equation}
   If further $u$ is nondegenerate in the sense that 
   $$
   \int_{Z_1} |\gd u|^2 \dd x \ge  \theta \int_{Z_2} |\gd u|^2 \dd x,
   $$
   then
   \begin{equation}
       \label{eq.local.coercivity}
       \int_{I_1} |\gd u|^2 \dd \sigma + \int_{Z_2} |F|^2 \dd x \ge C_\theta \int_{Z_2}|\gd u|^2 \dd x.
   \end{equation}
   The constant $C_\theta$ depends at most on $d$, $\theta$, $\norm{A}_{C^1(Z_2)}$ and geometry of $\Omega$.
\end{lemma}

\begin{proof}
    The proof is an immediate adaptation of the compactness argument from \cite[Lemma 5.3]{kenig_estimates_2013} with $E = 0$. Note that the only new ingredient is the source term $F$, which in the compactness argument vanishes in the limit equation, and hence the same unique continuation principle works. We also remark that although the domain $Z_2$ we choose are half-ball shaped, which is different from the cylinders in \cite[Lemma 5.3]{kenig_estimates_2013}, the proof does not need much change.
\end{proof}

\begin{lemma}
    \label{l.lower.bound.zep}
    Assume the same hypotheses as in Lemma \ref{l.boundary.estimate.Z.ep}. Then for some constants $C,c>0$,
    \begin{equation}
        \label{eq.lower.bound.zep}
         \int_{\pt\Om}|\gd z_\ep|^2 \dd \sigma \ge   c \lambda - C\ep \int_{\Om_{4\eps}}|\La_\eps z_\eps|^2 \dd x.
    \end{equation}
\end{lemma}

\begin{proof}
    We choose a maximally $\eps$-separated set $S_\ep\subset \pt\Om$, i.e. every two distinct points $p,q\in S_\eps$ satisfy $\dist(p,q)\ge \eps$, and adding one point $e\in \pt \Om \setminus S_\eps$ into $S_\eps$ breaks this property. To each $p\in S_\eps$ we associate  balls $D_s(p) : = B_{s \eps}(p)$. Note that the balls $D_s(p)$ satisfies: 
    \begin{enumerate}
        \item  \label{item.1.43} $\{D_2(p)\}_{p\in S_\eps}$ covers $\pt\Om$;
        
        \item  \label{item.2.43} $\bigcup_{p\in S_\eps} D_4(p)\cap \Om \subset \Om_{4\eps}$;
        
        \item  \label{item.3.43} $\Om_{\eps}\subset \bigcup_{p\in S_\eps} D_2(p) \cap \Omega$;
        
        \item \label{item.4.43} there are constants $\ep_0,C>0$ depending only on the geometry of $\Om$ and dimension $d$ such that for all $\ep<\ep_0$
        \begin{equation}
            \label{eq.upperbound.multiplicaity}
            \sum_{p\in S_\eps}\one_{D_4(p)} (x) \le C \textup{ for all }x\in \Om.
        \end{equation}
    \end{enumerate}

    We call $D_2(p)$ a $\theta$-good ball if
    \begin{equation}
        \label{eq.theta.good.ball}
        \int_{D_2(p)\cap \Om} |\gd z_\ep|^2 \dd x \ge \theta \int_{D_4(p) \cap \Om} |\gd z_\ep|^2 \dd x.
    \end{equation}
    Note that by Lemma \ref{l.boundary.estimate.Z.ep} and items (\ref{item.2.43}) - (\ref{item.4.43}) above, there is 
    \begin{equation}
        \label{eq.sum.control.good.ball}
        \sum_{p\in S_\eps} \frac{1}{\eps}\int_{D_2(p)\cap \Om} |\gd z_\ep|^2 \dd x \ge c \lambda \quad \textup{and} \quad \sum_{p\in S_\eps} \frac{1}{\eps}\int_{D_4(p)\cap \Om} |\gd z_\ep|^2 \dd x \le C \lambda.
    \end{equation}
On the other hand, the energy on the $\theta$-bad balls satisfies
\begin{equation}
    \label{eq.bad.balls.energy}
    \sum_{p\in S_\eps:\theta\textup{-bad}}\frac{1}{\eps} \int_{D_2(p)\cap \Om} |\gd z_\ep|^2 \dd x  \le \theta \sum_{p\in S_\eps} \frac{1}{\eps}\int_{D_4(p)\cap \Om} |\gd z_\ep|^2 \dd x \le C \theta\lambda.
\end{equation}
This shows that by choosing $\theta= \frac{c}{2C}$
\begin{equation}
    \label{eq.good.bad}
     \sum_{p\in S_\eps: \theta\textup{-good}}\frac{1}{\eps} \int_{D_2(p)\cap \Om} |\gd z_\ep|^2 \dd x \ge \frac{c}{2} \lambda.
\end{equation}

By rescaling the estimates in Lemma \ref{l.local.coercivity} for each $\theta$-good ball $D_2(p)$, \eqref{eq.theta.good.ball} implies
    \begin{equation}
        \label{eq.rescaled.local.coercive}
        \int_{D_2(p) \cap \pt\Om} |\gd z_\eps|^2 \dd \sigma + \eps \int_{D_4(p) \cap \Om} |\La_\eps z_\eps|^2 \dd x \ge \frac{C_\theta}{\eps}\int_{D_4(p) \cap \Om} |\gd z_\eps|^2 \dd x.
    \end{equation}
By summing up \eqref{eq.rescaled.local.coercive} in $p\in S_\eps$ that are $\theta$-good, \eqref{eq.good.bad},  item (\ref{item.1.43}) and item (\ref{item.4.43}), we finish the proof of the lemma.
\end{proof}

Let us now finish the proof of Theorem \ref{intro:main} and Theorem \ref{intro:cluster}.

\begin{proof}[Proof of Theorem \ref{intro:main} and Theorem \ref{intro:cluster}]
By Theorem \ref{t.SimpliesU}, it suffices to prove Theorem \ref{intro:cluster}.

We observe that by Lemma \ref{l.boundary.collar.gradient.estimate.apx}, for all unit vector $z_\ep \in X_{\sqrt{\lambda},\rho}$ 
    \[
    \begin{split}
        \|\La_\eps z_\eps\|_{L^2(\Om_{4\ep})} &\le \lambda \norm{z_\eps}_{L^2(\Om_{4\eps})} + \norm{r_\eps}_{2}\\
        &\le C \eps^{1/2}\lambda \norm{z_\eps}_{2}^{1/2} \norm{z_\eps}_{H^1(\Om)}^{1/2} + C \rho \sqrt{\lambda}\\
        &\le C\eps^{1/2}\lambda^{1+1/4} + C \rho\sqrt{\lambda},
    \end{split}
      \]
      where we have used \eqref{eq.homog.v.43} and the bound for $r_\eps$ in \eqref{eq.bound.reps.43}. This shows that 
      \[
      \eps \int_{\Om_{4\ep}} |\La_\ep z_\eps|^2 \dd x \le C \lmb\ep^2 \lambda^{3/2} + \eps\rho^2\rmb \cdot \lambda \le C (\delta^{3/2} + \eps \rho^2)\lambda.
      \]
      This finishes the proof by combining Lemma \ref{l.bl.blow.up} and Lemma \ref{l.lower.bound.zep}, and choosing $\delta, \eps_0, \rho_0$ sufficiently small.
\end{proof}

\section{Sharp commutator estimates}
\label{sec.sharp.commutator.estimate}

Let $\Omega$ be a bounded $C^{2,1}$ domain. In this section we consider the Dirichlet problem
\begin{equation}
    \cL_\e (u_\e) = 0 \quad \text{in } \Omega, \qquad u_\e = f \quad \text{on } \partial \Omega.
\end{equation}
Denote by $\cH_\e: f \to u_\e$ the map from boundary data $f \in L^2(\partial \Omega)$ to the solution $u_\e$ given by the Poisson integral (i.e., $\Lop_\eps$-harmonic extension). That is,
\begin{equation}\label{eq:Poisson-integral}
    \cH_\e(f)(x) = \int_{\partial \Omega} P_\e(x,Y) f(Y) \dd \sigma(Y).
\end{equation}
Moreover, $\cH_0$ denotes the solution map given by the homogenized operator $\cL_0$. Up to an affine change of variable, we can assume $\cL_0 = -\Delta$.

Let $b \in L^\infty(\partial \Omega)$ and $\phi \in H^1(\Omega)$. Define the commutator
\begin{equation}\label{eq:def-commutator}
    [\cH_\e, \phi] b = \cH_\e(\phi b ) - \phi \cH_\e (b).
\end{equation}

\begin{theorem}\label{thm:commutator}
For every fixed $B>0$ there are $\e_B>0$ and $C_B<\infty$ such that,
for $0<\e<\e_B$, $b\in L^\infty(\partial\Omega)$ and $\phi \in H^1(\Omega)$,
the function $W_\e = [\cH_\e, \phi] b$ belongs to $H^1_0(\Omega)$ and satisfies
\begin{equation}\label{eq:main}
    \| \nabla W_\e \|_{L^2(\Omega_{B\e})}
    \le C_B \| b \|_{L^\infty(\partial \Omega)} \big( \e^{1/2} \| \phi \|_{L^2(\partial \Omega)} + \| \nabla \phi \|_{L^2(\Omega_{2B\e})} \big)
\end{equation}
The constants are independent of $b,\phi,\e$.
\end{theorem}

This estimate can be viewed as a large-scale commutator estimate in periodic homogenization setting, which is of independent interest. As $\e \to 0$, it assembles a Dirichlet-to-Neumann map commutator estimate.

Recall the expansion of Poisson kernel in \cite[Theorem 3.11]{KLS14}, for $x\in \Omega, Y\in \partial \Omega$,
\begin{equation}\label{eq:Poisson-Expand}
    \nabla_x P_\e(x,Y) = \nabla_x \Phi_\e (x)\nabla_x P_0(x,Y) \omega_\e(Y) + R_\e(x,Y),
\end{equation}
where
    \[
 \omega_\e(Y)
 =\frac{\nu(Y)\cdot A^\e(Y)\nu(Y)  }
             {\nu(Y)\cdot\widehat A \nu(Y)} ~\lb \nu(Y)\cdot\nabla \Phi_\e(Y) \nu(Y)\rb,
\]
and the remainder satisfies, with $r = |x-Y|$
\begin{equation}\label{eq:Rsize}
    |R_\eps(x,Y)|
    \le C\min\left\{r^{-d},\,
                  \eps\log(2+r/\eps)r^{-d-1}\right\}.
\end{equation}
Recall that $\omega_\e(Y)$ is a uniformly bounded function on $\partial \Omega$ and $\nabla \Phi_\e(x)$ is uniformly bounded in $\Omega$, independent of $\e$. The expansion \eqref{eq:Poisson-Expand} was proved in $C^{3,\alpha}$ domains in \cite{KLS14}. A careful check on their proof shows that the entire proof works for $C^{2,\alpha}$ domains, except $C^\alpha$ continuity of $\nabla_x^2 \nabla_y G_0(x,y)$ in $x$ was used in the proof. However, $C^{2,\alpha}$ domains are also sufficient for these mixed third derivatives of the homogenized Green's function. Actually, if we set $F(x,y) = \nabla_y G_0(x,y)$, then $\cL_0 (F(\cdot, y)) =0 $ in $\Omega \setminus \{ y \}$ and $F(x,y) = 0$ for $x \in \partial \Omega$. This is because taking derivatives in $y$ does not change the zero Dirichlet boundary condition. Therefore, the boundary $C^{2,\alpha}$ estimate in $C^{2,\alpha}$ domains implies that $\nabla_x^2 \nabla_y G_0(x,y)$ is $C^\alpha$ in $x$. In the present setting, $C^{2,1}$ is sufficient for \eqref{eq:Poisson-Expand} and the remainder estimate.

Using \eqref{eq:Poisson-integral}, \eqref{eq:def-commutator} and \eqref{eq:Poisson-Expand}, we can write
\begin{equation}
    \begin{aligned}\label{eq:DxW-reduction-1}
        \nabla_x W_\e(x) & = \int_{\partial \Omega} \nabla_x P_\e(x, Y) (\phi(Y) - \phi(x)) b(Y) \dd \sigma(Y) - \nabla \phi(x) \cH_\e(b) \\
        & = \nabla \Phi_\e(x) \int_{\partial \Omega} \nabla_x P_0(x, Y) \omega_\e(Y) (\phi(Y) - \phi(x))  b(Y) \dd \sigma(Y) \\
        & \qquad + \int_{\partial \Omega} R_\e(x,Y) (\phi(Y) - \phi(x))  b(Y) \dd \sigma(Y) - \nabla \phi(x) \cH_\e(b).
    \end{aligned}
\end{equation}
We need to establish the boundary layer estimates for all three terms. It turns out that the last two terms can be handled by \eqref{eq:Rsize} and the maximum principle, respectively. For the leading term, we note that $\nabla \Phi_\e(x)$ is uniformly bounded, and $\omega_\e(Y)$ can be combined with $b(Y)$ and viewed as a new $L^\infty$ function. Consequently, the problem is reduced to studying
\begin{equation}\label{eq:DxW-reduction-2}
    T(\phi,b)(x) := \int_{\partial \Omega} \nabla P_0(x,Y) (\phi(Y) - \phi(x)) b(Y) \dd \sigma(Y).
\end{equation}
Let $Z(x) = [\cH_0, \phi]b(x)$. Then it is easy to see that
\begin{equation}\label{eq:DxW-reduction-3}
    T(\phi,b)(x) = \nabla Z(x) + \nabla \phi(x) \cH_0 b.
\end{equation}
Again, the second term is easy to handle by the maximum principle. Therefore, it suffices to consider $\nabla Z$.

\begin{lemma}\label{lem:localenergy}
Let $Z = [\cH_0, \phi]b$. For $2t< \text{diam}(\Omega)$,
\begin{equation}\label{eq:localenergy}
 \int_{\Omega_t}|\nabla Z|^2
 \le C \| b \|_{L^\infty(\partial \Omega)}^2\int_{\Omega_{2t}}|\nabla \phi|^2
       +C t^{-2}\int_{\Omega_{2t}\setminus\Omega_t}|Z|^2.
\end{equation}
Also,
\begin{equation}\label{eq:globalenergy}
    \| \nabla Z \|_{L^2(\Omega)} \le C\| b \|_{L^\infty(\partial \Omega)} \| \nabla \phi \|_{L^2(\Omega)}.
\end{equation}
These statements apply to $[\cH_\e, \phi]b$ with constants independent of $\eps$.
\end{lemma}
\begin{proof}
Set $U = \cH_0(b)$. The maximum principle implies $\| U \|_{L^\infty(\Omega)} \le \|b \|_{L^\infty(\partial \Omega)}$. We first show that for any $\psi \in H^1_0(\Omega)$,
\begin{equation}\label{eq:Cacc-H0b}
    \int_{\Omega} |\nabla U|^2 |\psi|^2 \le C\| b \|_{L^\infty(\partial \Omega)}^2 \int_{\Omega} |\nabla \psi|^2.
\end{equation}
Actually, testing the equation $-\Delta U = 0$ by $U \psi^2$ with $\psi \in C_0^\infty(\Omega)$ and applying an integration by parts, we can get
\begin{equation*}
    \int_{\Omega} |\nabla U|^2 |\psi|^2 \le C \int_{\Omega} |U|^2 |\nabla \psi|^2.
\end{equation*}
The $L^\infty$ bound of $U$ and a density argument yield \eqref{eq:Cacc-H0b}.

Next, a direct computation gives
\begin{equation}\label{eq:Zpde}
 -\Delta Z = \diver( U \nabla \phi) + \nabla U \cdot\nabla \phi,
 \qquad Z|_{\partial\Omega}=0.
\end{equation}
Choose $\zeta=1$ on $\Omega_t$, $\zeta=0$ outside $\Omega_{2t}$ and
$|\nabla\zeta|\le C/t$. Testing \eqref{eq:Zpde} with $\zeta^2 Z$ gives
\begin{equation*}
    \int_{\Omega} |\nabla Z|^2 |\zeta|^2 = -2\int_{\Omega} \nabla Z \cdot Z \zeta \nabla \zeta - \int_{\Omega} U\nabla \phi \cdot \nabla (Z \zeta^2) + \int_{\Omega} \nabla U \cdot \nabla \phi Z \zeta^2.
\end{equation*}
The last integral with $\nabla U$ can be handled using \eqref{eq:Cacc-H0b},
\begin{equation*}
\begin{aligned}
    \bigg| \int_{\Omega} \nabla U \cdot \nabla \phi Z \zeta^2 \bigg| & \le \| \zeta \nabla \phi \|_{L^2(\Omega)} \| \nabla U (Z\zeta) \|_{L^2(\Omega)} \\
    & \le C \| \zeta \nabla \phi \|_{L^2(\Omega)} \| b \|_{L^\infty(\partial \Omega)} \| \nabla (Z\zeta) \|_{L^2(\Omega)}.
\end{aligned}
\end{equation*}
Hence, by a standard argument we obtain
\[
 \int_\Omega \zeta^2|\nabla Z|^2
 \le C \| b\|_{L^\infty(\partial \Omega)}^2 \int_\Omega \zeta^2|\nabla \phi|^2+C\int_\Omega |\nabla\zeta|^2|Z|^2,
\]
which proves \eqref{eq:localenergy} by our choice of $\zeta$. Taking $\zeta=1$ proves \eqref{eq:globalenergy}.
\end{proof}

We next recall a commutator estimate for Dirichlet-to-Neumann map. Let $D$ be a bounded Lipschitz domain and $\cH_0$ be the harmonic extension in $D$ as before. The Dirichlet-to-Neumann map is given by 
$\Lambda_D f = \frac{\partial}{\partial \nu} \cH_0 f$, initially defined for smooth functions on $\partial D$. It has been proved in \cite[Theorem~5.1]{KLS14} (see also \cite[Theorem~1.2]{Shen16})
\begin{equation}\label{eq:DNgoal}
    \| [\Lambda_D, \phi] b \|_{L^2(\partial D)} \le C\| b \|_{L^\infty(\partial D)} \| \phi \|_{H^1(\partial D)},
\end{equation}
where $[\Lambda_D,\phi]b=\Lambda_D(\phi b)-\phi\Lambda_D b$.

The commutator estimate of Dirichlet-to-Neumann map is the key ingredient for the following boundary layer estimate.
\begin{lemma}\label{lem:Z}
    Let $Z$ be as above. Then there exists a constant $s_*$ depending only on $\Omega$ such that for $t<s_*$,
    \begin{equation}
        \int_{\Omega_t} |Z|^2 \le Ct^2 \| b\|_{L^\infty(\partial \Omega)}^2 \int_{\Omega_t} (|\nabla \phi|^2 + |\phi|^2).
    \end{equation}
\end{lemma}

\begin{proof}
    Let $\Omega^s = \{ x\in \Omega: \dist(x, \partial \Omega) > s \}$.
    For small $s$, define the following map
    \begin{equation}\label{eq:Def-Xs}
        X_s(Q)=Q-s \nu(Q),\qquad Q\in\partial\Omega,
    \end{equation}
    where $\nu(Q)$ is the unit outer normal at $Q$ to $\partial \Omega$. Since $\Omega$ is a bounded $C^{2,1}$ domain, there exists $s_*>0$ such that for any $s\in (0, s_*)$, the map $X_s: \partial \Omega \mapsto \partial \Omega_s$ is a bi-$C^{1,1}$ diffeomorphism, and $X_s(\partial \Omega) = \partial \Omega^s$ is also a $C^{2,1}$ domain. The unit outer normal at $X_s(Q)$ to the surface $\partial \Omega^s$, write as $\nu_s(X_s(Q))$, is equal to $\nu(Q)$ for any $0<s<s_*$. Denote by $J_s$ the Jacobian of $X_s$, i.e., $\dd\sigma_s(X_s(Q)) = J_s \dd\sigma(Q)$, where $d\sigma_s$ be the surface measure on $\partial \Omega^s$. Clearly, $J_s \simeq 1$ and $|\dd J_s/\dd s| \le C$.

    Let $F \in C^1(\Omega)$ and $0<s<s_*$. We introduce a surface transport formula
    \begin{equation}\label{eq:ds-domain}
        \frac{\dd}{\dd s} \int_{\partial \Omega^s} F \dd \sigma_s = \int_{\partial \Omega^s} \bigg( - \frac{\partial}{\partial \nu_s} F +  \mathcal{J}_s F \bigg) \dd \sigma_s,
    \end{equation}
    where $\nu_s$ is the outward normal on $\partial \Omega^s$ and $\mathcal{J}_s(X_s(Q)) = \frac{\dd J_s}{\dd s} J_s^{-1}(Q)$. Indeed, write
    \begin{equation}
        \int_{\partial \Omega^s} F \dd \sigma_s = \int_{\partial \Omega} F(X_s(Q)) J_s(Q) \dd \sigma(Q).
    \end{equation}
    Taking derivative in $s$, noting $\dd X_s/\dd s = -\nu(Q) = -\nu_s(X_s(Q))$, and changing variable back to $\partial \Omega^s$, we obtain \eqref{eq:ds-domain}.

    Recall that $Z = \cH_0(\phi b) - \phi \cH_0(b)$. Let $U = \cH_0( b)$ and  $ V = \cH_0(\phi b)$. Thus both $U$ and $V$ are harmonic in $\Omega$. Then $Z = V - \phi U$ and $Z = 0$ on $\partial \Omega$. 
    Let
    \begin{equation}
        E(s) = \int_{\partial \Omega^s} Z^2 \dd \sigma_s.
    \end{equation}
    Applying \eqref{eq:ds-domain}, 
    \begin{equation}\label{eq:E'}
        E'(s) = \int_{\partial \Omega^s} \bigg( -2 Z \frac{\partial Z}{\partial \nu_s} + \mathcal{J}_s Z^2 \bigg) \dd \sigma_s.
    \end{equation}
    
    Let $\Lambda_s$ be the Dirichlet-to-Neumann map in $\Omega^s$. Since $U, V$ are harmonic functions in $\Omega^s$, we have $\Lambda_s U = \frac{\partial U}{\partial \nu_s}$ and $\Lambda_s V = \frac{\partial V}{\partial \nu_s}$. Hence, on $\partial \Omega^s$,
    \begin{equation}
    \begin{aligned}
        \frac{\partial Z}{\partial \nu_s} & = \frac{\partial V}{\partial \nu_s} - \phi \frac{\partial U}{\partial \nu_s} - \frac{\partial \phi}{\partial \nu_s} U \\
        & = \Lambda_s V - \phi \Lambda_s U - \frac{\partial \phi}{\partial \nu_s} U \\
        & = \Lambda_s (V - \phi U) + ( \Lambda_s(\phi U) - \phi \Lambda_s U) - \frac{\partial \phi}{\partial \nu_s} U \\
        & = \Lambda_s Z + [\Lambda_s, \phi] U - \frac{\partial \phi}{\partial \nu_s} U
    \end{aligned}
    \end{equation}
    Inserting the last equation into \eqref{eq:E'}, we get
    \begin{equation}
    \begin{aligned}\label{eq.Es'}
        E'(s) + 2\int_{\partial \Omega^s} Z \Lambda_s Z \dd \sigma_s & = -2\int_{\partial \Omega^s} Z [\Lambda_s, \phi] U \dd \sigma_s + 2\int_{\partial \Omega^s} \frac{\partial \phi}{\partial \nu_s} U Z \dd \sigma_s \\
        & \qquad + \int_{\partial \Omega^s}  \mathcal{J}_s Z^2 \dd \sigma_s.
        \end{aligned}
    \end{equation}
    The integral on the left-hand side is nonnegative, because, with $Z_s = Z|_{\partial \Omega^s}$
    \begin{equation*}
        \int_{\partial \Omega^s} Z \Lambda_s Z \dd \sigma_s = \int_{\partial \Omega^s} \cH_0^s Z_s \frac{\partial}{\partial \nu_s} \cH_0^s Z_s \dd \sigma_s = \int_{\Omega^s} |\nabla \cH_0^s Z_s|^2 \ge 0,
    \end{equation*}
    where $\cH_0^s$ is the harmonic extension in $\Omega^s$.
    For the commutator of Dirichlet-to-Neumann map, it follows from \eqref{eq:DNgoal} that
    \begin{equation}
        \| [\Lambda_s, \phi] U \|_{L^2(\partial \Omega^s)} \le C\| U\|_{L^\infty(\partial \Omega^s)} \| \phi \|_{H^1(\partial \Omega^s)}.
    \end{equation}
    Put $M = \| b\|_{L^\infty(\partial \Omega)} = \| U \|_{L^\infty(\Omega)}$ (maximum principle) and 
    $$Q(s)^2 = \int_{\partial \Omega^s} (|\phi|^2 + |\nabla \phi|^2) \dd \sigma_s.$$
    As a consequence of \eqref{eq.Es'},
    \begin{equation}
        E'(s) \le CM Q(s) E(s)^{1/2} + CE(s).
    \end{equation}
    Let $G(s) = E(s)^{1/2}$. Then
    \begin{equation}
        G'(s) \le CMQ(s) + C G(s).
    \end{equation}
    Hence, by $G(0) = 0$, for $s\in (0,s_*)$,
    \begin{equation}
        G(s) \le CM e^{Cs} \int_0^s Q(r) dr \le CM s^{1/2} \bigg( \int_0^s Q(r)^2 dr\bigg)^{1/2}.
    \end{equation}
    This, together with the coarea formula, gives
    \begin{equation}
        \int_{\partial \Omega^s} Z^2 \dd \sigma_s \le CM^2 s \int_{\Omega_s} (|\nabla \phi|^2 + |\phi|^2).
    \end{equation}
    Integrating $s$ over $(0,t)$ and by the coarea formula again,
    \begin{equation}
        \int_{\Omega_t} |Z|^2 = \int_0^t \int_{\partial \Omega^s} |Z|^2 \dd \sigma_s ds \le CM^2 t^2 \int_{\Omega_t} (|\nabla \phi|^2 + |\phi|^2).
    \end{equation}
    This ends the proof.
\end{proof}

\begin{proof}[Proof of Theorem \ref{thm:commutator}]
    We first recall the reductions from \eqref{eq:DxW-reduction-1} -- \eqref{eq:DxW-reduction-3}. Combining Lemma \ref{lem:localenergy} and Lemma \ref{lem:Z}, we have
    \begin{equation}
        \int_{\Omega_t} |\nabla Z|^2 \le C \| b\|_{L^\infty(\partial \Omega)}^2 \int_{\Omega_{2t}} (|\nabla \phi|^2 + |\phi|^2),
    \end{equation}
    for any $t<s_*/2$. The term $\nabla \phi(x) \cH_0 b$ has the same bound by the maximum principle. Thus
    \begin{equation}\label{eq:T-layer-est}
    \begin{aligned}
        \int_{\Omega_t} |T(\phi, b)|^2 & \le C \| b\|_{L^\infty(\partial \Omega)}^2 \int_{\Omega_{2t}} (|\nabla \phi|^2 + |\phi|^2) \\
        & \le C \| b\|_{L^\infty(\partial \Omega)}^2 \bigg( t \int_{\partial \Omega} |\phi|^2 \dd \sigma +  \int_{\Omega_{2t}} |\nabla \phi|^2 \bigg),
    \end{aligned}
    \end{equation}
    where in the last inequality we have used Lemma \ref{lem.layer2bdry}.
    Now note that the leading integral in \eqref{eq:DxW-reduction-1} is exactly $\nabla \Phi_\e T(\phi, \omega_\e b)$, which in view of \eqref{eq:T-layer-est} gives the desired bound with $t = B\e$ since $\nabla \Phi_\e$ and $\omega_\e$ are uniformly bounded.

    The last term $\nabla \phi(x) \cH_\e b$ in \eqref{eq:DxW-reduction-1} can also be handled by the maximum principle giving the same bound. It remains to estimate the second integral in \eqref{eq:DxW-reduction-1} with remainder $R_\e$.

    Let
    \begin{equation}
        K(x) = \int_{\partial \Omega} R_\e(x,Y) (\phi(Y) - \phi(x))  b(Y) \dd \sigma(Y).
    \end{equation}
    Put $M = \| b\|_{L^\infty(\partial \Omega)}$.
    Consider $x = X_s(Q)$ on $\partial \Omega^s$.
    We split the integral into two parts $K(x) = K_{\rm nor}(x) + K_{\tan}(x)$, where
    \begin{align}
        K_{\rm nor}(X_s(Q)) & = \int_{\partial \Omega} R_\e(X_s(Q),Y) (\phi(Q) - \phi(X_s(Q)))  b(Y) \dd \sigma(Y), \\
        K_{\tan}(X_s(Q)) &= \int_{\partial \Omega} R_\e(X_s(Q),Y) (\phi(Y) - \phi(Q))  b(Y) \dd \sigma(Y).
    \end{align}

    We first estimate the normal term $K_{\rm nor}(X_s(Q))$. By the first estimate of \eqref{eq:Rsize}, we have
    \begin{equation}
        \sup_{Q\in \partial \Omega} \int_{\partial \Omega} |R_\e(X_s(Q), Y)| \dd \sigma(Y) + \sup_{Y\in \partial \Omega} \int_{\partial \Omega} |R_\e(X_s(Q), Y)| \dd \sigma(Q) \le \frac{C}{s}.
    \end{equation}
    By these and Schur test, we have
    \begin{equation}
    \begin{aligned}
        \| K_{\rm nor}(X_s(\cdot)) \|_{L^2(\partial \Omega)} &\le CM s^{-1} \| \phi(\cdot) - \phi(X_s(\cdot)) \|_{L^2(\partial \Omega)} \\
        & \le CM s^{-1} \int_0^s \| \nabla \phi(X_t(\cdot)) \|_{L^2(\partial \Omega)} \dd t,
    \end{aligned}
    \end{equation}
    where we also used
    $$
    \phi(X_s(Q)) - \phi(Q) = \int_0^s -n\cdot \nabla \phi (X_t(Q)) dt
    $$
    in the last inequality.
    Then, by coarea formula and Hardy's inequality (see \eqref{eq:Hardy-1} below), we arrive at
    \begin{equation}
    \begin{aligned}
        \int_{\Omega_{T}} |K_{\rm nor}|^2 & \simeq \int_0^T \| K_{\rm nor}(X_s(\cdot)) \|_{L^2(\partial \Omega)}^2 \dd s \\
        & \le CM^2 \int_0^T \| \nabla \phi(X_t(\cdot)) \|_{L^2(\partial \Omega)}^2 \dd t \simeq CM^2 \int_{\Omega_T} |\nabla \phi|^2,
    \end{aligned}
    \end{equation}
    for $T < s_*$.

    Next, we estimate tangential term $K_{\tan}(X_s(Q))$. 
    
    We need a simple fact first given in flat $\R^{d-1}$. Let $r > 0$. Suppose that $f \in C^1 \cap H^1(\R^{d-1})$. Then there exists $C>0$ depending only on $d$ such that
    \begin{equation}
        \| \fint_{B_r(Q)} |f(Y) - f(Q)| \dd Y \|_{L^2(\R^{d-1}; \dd Q)} \le C r \| \nabla f \|_{L^2(\R^{d-1})}.
    \end{equation}
    This property can be extended to $\partial \Omega$ by a partition of unity, i.e., for $r < r_0$ for some $r_0$ depending only on $\Omega$, and $f\in H^1(\partial \Omega)$, we have
    \begin{equation}\label{eq:surface-LS}
        \| \fint_{\partial \Omega \cap B_r(Q)} |f(Y) - f(Q)| \dd \sigma(Y) \|_{L^2(\partial \Omega; \dd \sigma(Q))} \le C r \| \nabla f \|_{L^2(\partial \Omega)}.
    \end{equation}
    Without loss of generality, assume $r_0 = s_*$.

    For $t\in (0,s_*)$, define
    \begin{equation}
        a(t) = \| \nabla \phi(X_t) \|_{L^2(\partial \Omega)} \simeq \| \nabla \phi \|_{L^2(\partial \Omega^t)}.
    \end{equation}
    Write
    \begin{equation}
    \begin{aligned}
        \phi(Y) - \phi(Q) & = \phi(X_t(Y)) - \phi(X_t(Q)) \\
        &\qquad -\int_0^t n\cdot \nabla \phi(X_\tau(Q)) \dd\tau + \int_0^t n\cdot \nabla \phi(X_\tau(Y)) \dd \tau 
    \end{aligned}
    \end{equation}
    Combining this and \eqref{eq:surface-LS}, we have
    \begin{equation}\label{eq:trace-1}
        \| \fint_{\partial \Omega \cap B_r(Q)} |\phi(Y) - \phi(Q)| \dd \sigma(Y) \|_{L^2(\partial \Omega)} \le Cr a(t) + C\int_0^t a(\tau) \dd \tau,
    \end{equation}
    for any $t\in (0, s_*)$.
    Taking average in $t$ over $(0,r)$, we then obtain
    \begin{equation}\label{eq:trace-average}
        \| \fint_{\partial \Omega \cap B_r(Q)} |\phi(Y) - \phi(Q)| \|_{L^2(\partial \Omega)} \le C \int_0^r a(\tau) \dd \tau.
    \end{equation}
    
    Let $T = 2B\e \le s_*$. Fix $s \in (0,T)$. Write
    \begin{equation}
    \begin{aligned}
        & K_{\tan}(X_s(Q))  \\
        & = \int_{\partial \Omega \cap \{ |Y-Q|< s \}}  R_\e(X_s(Q),Y) (\phi(Y) - \phi(Q))  b(Y) \dd \sigma(Y) \\
        & + \sum_{j \ge 0} \int_{\partial \Omega \cap \{ 2^js \le |Y-Q|<2^{j+1}s \}}  R_\e(X_s(Q),Y) (\phi(Y) - \phi(Q))  b(Y) \dd \sigma(Y).
    \end{aligned}
    \end{equation}
    Observe that $|X_s(Q) - Y| \simeq s + |Q - Y|$.
    Now, we divide $K_{\tan}$ into three parts: $2^j s \le B\e,\  B\e < 2^j s \le r_0$ and $2^j s > r_0$, denoted by $K^{<B\e}_{\tan}$, $K^{(B\e,r_0)}_{\tan}$ and $K^{>r_0}_{\tan}$, respectively. If $s > B\e$, then $K^{<B\e}_{\tan} = 0$. Otherwise,
    \begin{equation}
    \begin{aligned}
        |K^{<B\e}_{\tan}(X_s(Q))|
         \le \sum_{2^j s \le B\e} \frac{CM}{2^j s} \fint_{\partial \Omega \cap B(Q,2^{j+1} s)} |\phi(Y) - \phi(Q)| \dd \sigma(Y).
    \end{aligned}
    \end{equation}
    Taking $L^2$ norm in $Q \in \partial \Omega$ and applying \eqref{eq:trace-average} and the first bound of $R_\e$ in \eqref{eq:Rsize} give
    \begin{equation}
    \begin{aligned}\label{eq:Ktan-1}
        \| K^{<B\e}_{\tan}(X_s(\cdot)) \|_{L^2(\partial \Omega)} & \le \sum_{2^j s \le B\e} \frac{CM}{2^j s} \int_{0}^{2^{j+1} s} a(\tau) \dd \tau \le CM \int_0^{2B\e} \frac{a(t)}{t+s} \dd t \\
        & \le CM \frac{1}{s}\int_0^s a(t) \dd t + CM \int_s^{2B \e} \frac{a(t)}{t} \dd t.
    \end{aligned}
    \end{equation}

    For the second part, using the second bound of $R_\e$ in \eqref{eq:Rsize} and \eqref{eq:trace-1} gives
    \begin{equation}
    \begin{aligned}\label{eq:Ktan-2}
        \| K_{\tan}^{(B\e,r_0)}(X_s(\cdot)) \|_{L^2(\partial \Omega)} &\le \sum_{2^j s > B\e} CM \frac{\e \log(2 + 2^{j+1}s/\e)}{(2^j s)^2} \Big(2^{j+1}s a(s) + \int_0^s a(t) \dd t \Big) \\
        & \le CM \Big( a(s) + \frac{1}{s} \int_0^s a(t) \dd t \Big).
    \end{aligned}
    \end{equation}

    For the last part, as $2^j s > r_0$ (thus $|Y - Q| > r_0$), we have $|R_\e(X_s(Q), Y)| \le C$. Hence,
    \begin{equation*}
        |K^{>r_0}_{\tan}(X_s(Q))| \le CM |\phi(Q)| + CM \int_{\partial \Omega} |\phi| \dd \sigma,
    \end{equation*}
    yields
    \begin{equation}\label{eq:Ktan-3}
        \| K^{>r_0}_{\tan}(X_s(\cdot)) \|_{L^2(\partial \Omega)} \le CM \| \phi \|_{L^2(\partial\Omega)}.
    \end{equation}

    Finally, by Hardy's inequality and its dual form, we have
    \begin{equation}\label{eq:Hardy-1}
        \int_0^T \Big( \frac{1}{s} \int_0^s a(t) \dd t \Big)^2 \dd s \le 4\int_0^T a(t)^2 \dd t,
    \end{equation}
    and
    \begin{equation}\label{eq:Hardy-2}
        \int_0^T \Big( \int_s^{T} \frac{a(t)}{t} \dd t \Big)^2 \dd s \le 4\int_0^T a(t)^2 \dd t.
    \end{equation}
    Combining \eqref{eq:Ktan-1} -- \eqref{eq:Hardy-2}, obtain
    \begin{equation}
        \int_{\Omega_T} |K|^2 \simeq \int_0^T \| K(X_s(\cdot)) \|_{L^2(\partial \Omega)}^2 \dd s \le CM^2 \int_{\Omega_T} |\nabla \phi|^2 + CM^2 T \int_{\partial \Omega} |\phi|^2 \dd \sigma.
    \end{equation}
    This ends the proof.
\end{proof}

\appendix

\section{Technical lemmas}

All the following technical lemmas are established in a general bounded Lipschitz domain $\Omega$ in $\R^d$. 

\begin{lemma}
    \label{l.multiplicative.trace.theorem.apx}
    There exists $C>0$ depending only on $\Omega$ such that for all $v\in H^1(\Omega)$,
    \begin{equation}
        \label{eq.multiplicative.trace.theorem}
        \norm{v}_{L^2(\pt\Omega)} \le C \lb\norm{v}_{L^2(\Omega)}\norm{v}_{H^1(\Omega)}\rb^{1/2}
    \end{equation}
\end{lemma}

\begin{proof}
  Since $\Omega$ is a bounded Lipschitz domain there is a smooth vector field $h$ on $\overline{\Omega}$ such that $h\cdot \nu \ge c_0 >0$ uniformly on the boundary $\pt \Omega$, where $\nu$ is the unit outer normal defined at almost every boundary point. Now we have for all smooth function $v \in C^\infty(\overline{\Omega})$ 
  \begin{equation}
      \label{eq.integrate.by.parts.technical}
      \begin{split}
              \int_{\pt \Omega} |v|^2 \dd \sigma & \le C \int_{\pt \Omega} (h,\nu)|v|^2 \dd  \sigma\\
              & = C \int_{\Omega} \diver (h|v|^2) \dd x \\
&\le C \int_{\Omega} |v|^2 + |v||Dv| \dd x,
      \end{split}
    \end{equation}
    which proves the lemma by Cauchy-Schwarz inequality and a standard density argument.
\end{proof}

\begin{lemma}
    \label{l.boundary.collar.gradient.estimate.apx}
    There exists $t_0 > 0$ depending only on the geometry of $\Omega$, such that for all $v \in H^1(\Omega)$ and $0<t<t_0$, 
    \begin{equation}
    \label{eq.boundary.collar.gradient.estimate.apx}
        \frac{1}{t}\int_{\Om_t} |v|^2 \dd x \le  C \norm{v}_{L^2(\Omega)}\norm{v}_{H^1(\Omega)}.
    \end{equation}
\end{lemma}
\begin{proof}
Recall that $\Omega^s = \{x\in \Omega: \dist(x,\partial \Omega) > s \}$ is a Lipschitz domain with uniform Lipschitz constant. By Lemma \ref{l.multiplicative.trace.theorem.apx},
\begin{equation}
    \int_{\partial \Omega^s} |v|^2 \dd \sigma_s \le C\| v \|_{L^2(\Om^s)} \| v \|_{H^1(\Om^s)} \le C\| v \|_{L^2(\Om)} \| v \|_{H^1(\Om)}.
\end{equation}
Integrating over $s\in (0,t)$ and by coarea formula, we obtain the desired estimate.

\end{proof}

\begin{lemma}\label{lem.layer2bdry}
    There exists $C>0, t_0>0$ depending only on $\Omega$ such that for all $v\in H^1(\Omega)$ and $0<t<t_0$,
    \begin{equation}
        \norm{v}_{L^2(\Omega_t)} \le C \lb t^{1/2} \norm{v}_{L^2(\partial \Omega)} + t \norm{\nabla v}_{L^2(\Omega_t)} \rb,
    \end{equation}
    and
    \begin{equation}
        \norm{v}_{L^2(\partial \Omega)} \le C\lb t^{-1/2} \norm{v}_{L^2(\Omega_t)} + t^{1/2} \norm{\nabla v}_{L^2(\Omega_t)} \rb.
    \end{equation}
\end{lemma}
\begin{proof}
    Since $\Omega$ is a Lipschitz domain, if $t < t_0$ for some small $t_0>0$, $\Omega^t$ is still a Lipschitz domain with uniform constant. Moreover, we can find a smooth vector field $h$ on $\overline{\Om}$ such that $c_0 \le h\cdot  \nu_t \le C_0$ on $\partial \Omega^t$ for every $t\in (0,t_0)$. Such $h$ can be first constructed first locally in a graph domain and then glued by a partition of unity. Fix $t \in (0,t_0)$, it follows from the divergence theorem as before, for each $s\in (0,t)$,
    \begin{equation}\label{eq.2bdry-Rellich}
        \begin{aligned}
            \int_{\partial \Omega^s} (h\cdot \nu_s) |v|^2 \dd\sigma_s - \int_{\partial \Om} (h\cdot \nu) |v|^2 \dd \sigma = -\int_{\Omega_s} \diver( h |v|^2 ),
        \end{aligned}
    \end{equation}
    where $\Omega^s$ is given as in the proof of Lemma \ref{l.boundary.collar.gradient.estimate.apx}.
    Thus, by $h\cdot \nu_s \ge c_0$ and the Cauchy-Schwarz inequality with arbitrarily $\delta>0$,
    \begin{equation*}
        \int_{\partial \Omega^s} |v|^2 \dd \sigma_s \le C\int_{\partial \Om} |v|^2 \dd \sigma + C(1 + \delta t^{-1}) \int_{\Omega_t} |v|^2 + C \delta^{-1} t \int_{\Omega_t} |\nabla v|^2.
    \end{equation*}
    Integrating over $s \in (0,t)$, we obtain by the coarea formula
    \begin{equation*}
        \int_{\Omega_t} |v|^2 \le Ct \int_{\partial \Om} |v|^2 \dd \sigma + C (t+\delta) \int_{\Omega_t} |v|^2 + C\delta^{-1} t^2 \int_{\Omega_t} |\nabla v|^2.
    \end{equation*}
    By letting $t_0$ and $\delta>0$ small such that $C(t_0+\delta) < 1/2$, the second term on the right-hand side can be absorbed to the left-hand side. This implies the first desired estimate. The second estimate also follows from \eqref{eq.2bdry-Rellich} by using $h\cdot \nu \ge c_0$.
\end{proof}

\begin{lemma}
    \label{l.H1half.boundary.trace.apx}
    For any fixed Lipschitz 1-periodic function $f$ and $g \in H^1(\Omega)$, the following trace estimate holds
    \begin{equation}
        \label{eq.H1half.boundary.trace}
        \norm{f^\eps g}_{H^{1/2}( \partial\Omega)} \le C \lb \ep^{-1/2} (\norm{g}_{L^2(\Omega)} \norm{g}_{H^{1}(\Omega)})^{1/2} + \norm{g}_{H^{1}(\Omega)} \rb,
    \end{equation}
    where the constant $C>0$ may depend on $\norm{f}_{W^{1,\infty}(\T^d)}$ but is independent of $g$ and $\ep$. In particular, if $B_\ep$ is the $\Lop_\ep$-harmonic extension of $f^\eps g|_{\partial \Om}$, then
    \begin{equation}\label{eq.Be-energy}
        \| B_\ep \|_{H^1(\Omega)} \le C \lb \ep^{-1/2} (\norm{g}_{L^2(\Omega)} \norm{g}_{H^{1}(\Omega)})^{1/2} + \norm{g}_{H^{1}(\Omega)} \rb.
    \end{equation}
\end{lemma}

\begin{proof}

Let $G$ be the harmonic extension of $g$. By the classical $L^2$ estimate of nontangential maximal function in Lipschitz domains \cite[Corollary 2.1.6]{kenig_harmonic_1994}, we have
\begin{equation}\label{eq.G-nontangential}
    \| (G)^* \|_{L^2(\partial \Omega)} \le C\| g\|_{L^2(\partial \Omega)}.
\end{equation}
Moreover, the trace theorem implies $g\in H^{1/2}(\partial \Omega)$. Thus the elliptic energy estimate yields
\begin{equation}\label{eq.G-energy}
    \| G \|_{H^1(\Omega)} \le C\| g \|_{H^{1/2}(\partial \Omega)} \le C\|g  \|_{H^1(\Omega)}.
\end{equation}

For $\eps$ sufficiently small, Let $\phi_\ep \in C_0^\infty(\Omega)$ be a cutoff function, satisfying  $0\le \phi_\ep \le 1$ and
\[
\phi_\ep \equiv 0 \textup{ in }\Omega_{\ep}, \ \phi_\ep \equiv 1 \textup{ in }\Omega\setminus\Omega_{2\ep} \ \textup{ and }\ |\nabla \phi_\ep| \le C /\ep.
\]
By the trace theorem
\begin{equation}
    \label{eq.calculate.H1halfnorm}
    \begin{split}
         \norm{f^\eps g}_{H^{1/2}(\partial\Omega)} = \norm{f^\eps G}_{H^{1/2}(\partial\Omega)} &\le C\norm{f^\ep G (1-\phi_\ep)}_{H^1(\Omega)} \\
& \le C\eps^{-1} \norm{G}_{L^2(\Omega_{2\ep})}  + C\norm{\nabla G}_{L^2(\Omega_{2\ep})}\\
&\le C\eps^{-1/2} \norm{(G)^*}_{L^2(\partial\Omega)} + C\| g \|_{H^{1/2}(\partial \Omega)} \\
&\le C\eps^{-1/2} \norm{g}_{L^2(\partial\Omega)} + C\| g \|_{H^1(\Omega)},
    \end{split}
   \end{equation}
where we have used \eqref{eq.G-energy} and \eqref{eq.G-nontangential} in the last inequalities. This proves \eqref{eq.H1half.boundary.trace} in view of Lemma \ref{l.multiplicative.trace.theorem.apx}. Finally, \eqref{eq.Be-energy} follows from the elliptic energy estimate.
\end{proof}

\begin{lemma}\label{l.fine.estimate.rellich.integral.apx}
    Let $G_\ep$ denote the $\Lop_\eps$-harmonic extension of $g$ in $\Omega$, $h$ any smooth vector field and $v\in H_0^1(\Om)$. Then the following estimate holds
    \begin{equation}
    \label{eq.fine.estimate.rellich.integral.apx}
        \lw \int_{\Om} G_\ep h\cdot \nabla v\dd x\rw \le C \norm{g}_{L^2(\pt\Om)}\lma\norm{v}_{L^2(\Om)}^{1/2} \norm{v}_{H^1(\Om)}^{1/2}+\norm{v}_{L^2(\Om)}\rma.
    \end{equation}
\end{lemma}

\begin{proof}
First of all, by the $L^2$-estimate for the nontangential maximal function for the operator $\Lop_\eps$ in Lipschitz domains \cite[Theorem~1.1]{KS11}, we have
\begin{equation}\label{eq.Ge-L2nontangential}
    \norm{ (G_\ep)^*}_{L^2(\partial \Omega)} \le C\norm{g}_{L^2(\partial \Omega)}.
\end{equation}
If $\norm{v}_{H^1(\Omega)} \le C_1 \norm{v}_{L^2(\Omega)}$ for some fixed constant $C_1 > 0$, then Cauchy-Schwarz inequality and \eqref{eq.Ge-L2nontangential} directly imply the desired estimate. Hence, below we only consider the case $\norm{v}_{L^2(\Omega)}/\norm{v}_{H^1(\Omega)} \le c$ for some small $c$ depending on $\Omega$.

Let $t>0$ to be determined.
We decompose the integral into $\Omega = \Omega_t \cup \lb \Omega \setminus \Om_t\rb$.
Then we have
\begin{equation}
    \label{eq.compute.integral}
    \begin{split}
        \lw \int_{\Om} G_\ep h\cdot \nabla v\dd x\rw & \le \lw \int_{\Om_t} G_\ep h\cdot \nabla v\dd x\rw + \lw \int_{\Om\setminus\Om_t} G_\ep h\cdot \nabla v\dd x\rw  \\
& = : T_t + S_t.
    \end{split}
    \end{equation}
Note that by Cauchy-Schwarz and \eqref{eq.Ge-L2nontangential}, we have
\[
|T_t| \le C \norm{v}_{H^1} \norm{ G_\e}_{L^2(\Om_t)} \le C t^{1/2} \norm{g}_{L^2(\pt \Om)} \norm{v}_{H^1(\Om)}.
\]

On the other hand, by the Caccioppoli inequality and \eqref{eq.Ge-L2nontangential}
\[
\int_{\Om\setminus\Om_t}|\nabla G_\ep|^2 \dd x \le  \frac{C}{t^2} \int_{\Omega_t} |G_\eps|^2 \le \frac{C}{t} \norm{g}_{L^2(\pt \Om)}^2.
\]
This shows that after integration by parts and using \eqref{eq.Ge-L2nontangential}
\[
\begin{split}
    |S_t| &\le  \lw\int_{\Om\setminus\Om_t} \nabla G_\ep \cdot h v \dd x \rw + \lw\int_{\pt(\Omega \setminus \Omega_t)} G_\ep v (h\cdot \nu_t) \dd \sigma_t \rw + \lw \int_{\Om\setminus\Om_t} G_\ep v ~\textup{div}(h) \dd x \rw\\
    &\le Ct^{-\frac{1}{2}}\norm{g}_{L^2(\pt\Om)} \norm{v}_{L^2(\Om)} +C \norm{ g}_{L^2(\pt \Om )} \norm{v}_{L^2(\pt(\Om \setminus \Omega_t))}\\
    &\quad+C\norm{g}_{L^2(\pt\Om)}\norm{v}_{L^2(\Om)}\\
    &\le  Ct^{-\frac{1}{2}}\norm{g}_{L^2(\pt\Om)} \norm{v}_{L^2(\Om)} +C \norm{ g}_{L^2(\pt\Om)} \lb\norm{v}_{L^2(\Om)}\norm{v}_{H^1(\Om)}\rb^{1/2}\\
    &\quad+C\norm{g}_{L^2(\pt\Om)}\norm{v}_{L^2(\Om)},
\end{split}
\]
where we have used Lemma \ref{l.multiplicative.trace.theorem.apx} in the last line. The proof of the lemma is finished by combining the estimates for $T_t$ and $S_t$ above and choosing $t =  \norm{v}_{L^2(\Om)}/\norm{v}_{H^1(\Om)}$.  
\end{proof}

\bibliographystyle{abbrv}
\bibliography{ref}

\end{document}